\documentclass[a4paper,11pt]{amsart}      

\usepackage{graphicx}
\usepackage{multirow,yfonts}
\usepackage{dsfont}
\usepackage[applemac]{inputenc}
\usepackage[english]{babel}
\usepackage{amssymb, amsfonts}
\usepackage{amsmath,amsthm}
\usepackage{stmaryrd}
\usepackage{xypic}
\usepackage{graphicx}
\usepackage{chngpage}
\usepackage[toc,page]{appendix} 
\usepackage{subfigure}
\usepackage{color}
\usepackage{comment}
\usepackage{hyperref}
\usepackage{amsmath}
\usepackage{accents}

\definecolor{gris}{gray}{0.45}
\DeclareMathOperator{\dlog}{dlog}

\newcommand{\gf}{\Lambda_n}

\newcommand {\M}{\mathrm{MCG}}
\newcommand {\PM}{\mathrm{PMCG}}

\newcommand{\F}{\mathcal{F}}

\newcommand{\C}{\mathbb{C}} 
\newcommand{\D}{\mathbb{D}} 
 
\newcommand{\Z}{\mathbb{Z}}

\newcommand{\Pu}{\mathbb{P}^1}

            \newcommand{\overbull}[1]{\accentset{\bullet}{#1}}
                      \newcommand{\rev}{f}

                      \subjclass[2010]{14F35; 37F75}
                      \keywords{Logarithmic flat connections; Fundamental group representations}
\begin{document}

\theoremstyle{plain}
\newtheorem{theorem}{Theorem}[subsection]
\newtheorem{proposition}[theorem]{Proposition}
\newtheorem{lemma}[theorem]{Lemma}
\newtheorem{corollary}[theorem]{Corollary}
\newtheorem{question}{Question}
\newtheorem{THM}{Theorem}
\renewcommand*{\theTHM}{\Alph{THM}}
\newtheorem{COR}[THM]{Corollary}{\bf}{\itshape}
\newtheorem*{THM:maintheorem}{Theorem A'}{\bf}{\itshape}  
\newtheorem*{COR:1}{Corollary B}{\bf}{\itshape}  
\newtheorem*{COR:2}{Theorem C}{\bf}{\itshape}  
\newtheorem*{lem:A1}{Lemma A1}{\bf}{\itshape}  
\newtheorem*{lem:A2}{Lemma A2}{\bf}{\itshape}  
\newtheorem*{proofbis}{Proof}{\itshape}{\upshape}

\theoremstyle{definition}   
\newtheorem{definition}[theorem]{Definition}                 
\newtheorem{cond}{Condition}          
\newtheorem{prop/def}[theorem]{Proposition/Definition}
\newtheorem{example}[theorem]{Example}

\theoremstyle{remark}    
\newtheorem{remark}[theorem]{Remark}     

\title[Algebraic isomonodromic deformations]{Algebraic isomonodromic deformations of logarithmic connections on the Riemann sphere and finite braid group orbits on character varieties}

\author{Ga\"el Cousin}



\date{}

\maketitle
\vspace{-1cm}
\begin{abstract}
We study algebraic isomonodromic deformations of flat logarithmic connections on the Riemann sphere with $n\geq 4$ poles, for arbitrary rank. We introduce a natural property of algebraizability for the germ of universal deformation of such a connection. We relate this property to a peculiarity of the corresponding monodromy representation: to yield a finite braid group orbit on the appropriate character variety. Under reasonable assumptions on the deformed connection, we may actually establish an equivalence between both properties. We apply this result in the rank two case to relate finite branching and algebraicity for solutions of Garnier systems.

For general rank, a byproduct of this work is a tool to produce regular flat meromorphic connections on vector bundles over projective varieties of high dimension.
\keywords{Logarithmic flat connections \and Fundamental group representations}
\end{abstract}
\tableofcontents
\section{Introduction}
In this paper, we study isomonodromic deformations of logarithmic connections with $n\geq 4$ poles on the Riemann sphere.
From Malgrange \cite{MR0728431}, it is known that every such deformation parametrized by a simply connected basis is obtained by pull-back of a \emph{universal isomonodromic deformation}.

We introduce a notion of algebraizability for the germ of universal isomonodromic deformation of a given connection $\nabla$. Roughly speaking, this germ is considered algebraizable if we can recover it by restriction of a flat logarithmic connection on a (necessarily algebraic, by GAGA) vector bundle over a projective manifold $X$ to a suitable family of rational curves in $X$.  

Our main result relates  this algebraizability to a dynamical property for the monodromy representation of the deformed connection $\nabla$.

\begin{THM}[Main theorem]\label{maintheorem} 
Let $\nabla$ be a logarithmic connection on a rank $m$ holomorphic vector bundle over $\Pu$ with $n$ poles $x_1,\ldots,x_n$, $n\geq 4$.  

Let $\rho : \pi_1(\Pu\setminus \{x_1,\ldots,x_n\},x_0)\rightarrow \mathrm{GL}_m(\C)$ be its monodromy representation.
Suppose the germs of $\nabla$ at the poles are mild transversal models and $\rho$ is semisimple or $m=2$.
Then the following are equivalent.
\begin{enumerate}
\item \label{unthmintro}The conjugacy class $[\rho]$ has a finite orbit under the pure mapping class group $\mathrm{P}\M_n(\Pu)$.
\item \label{deuxthmintro} The germ of universal isomonodromic deformation of $\nabla$ is algebraizable.
\end{enumerate}
\end{THM}
The notion of mild transversal model is introduced in section $\ref{secRH}$.
Actually, if $\nabla$ is constructed from $\rho$ \textit{via} Deligne's canonical extension, the mildness assumption of Theorem \ref{maintheorem} is automatically satisfied.

If we work up to birational gauge transformations, we can drop this assumption.
\begin{COR}\label{COR:1}
 Let $\nabla$ be a logarithmic connection on a rank $m$ holomorphic vector bundle over $\Pu$ with $n$ poles $x_1,\ldots,x_n$, $n\geq 4$.  

Let $\rho : \pi_1(\Pu\setminus \{x_1,\ldots,x_n\},x_0)\rightarrow \mathrm{GL}_m(\C)$ be its monodromy representation.
Suppose  $\rho$ is semisimple or $m=2$.
Then the following are equivalent.
\begin{enumerate}
\item The conjugacy class $[\rho]$ has a finite orbit under the pure mapping class group $\mathrm{P}\M_n(\Pu)$.
\item   Up to a birational gauge transformation of $\nabla$, the germ of universal isomonodromic deformation of $\nabla$ is algebraizable.
\end{enumerate}
 \end{COR}
The introduction of technicalities concerning mild tranversal models is essential in the proof of Theorem \ref{COR:2} below, which is a generalization of the global result of Iwasaki \cite{MR2388081} concerning algebraic solutions of Painlev\'e VI equation.
\begin{THM}\label{COR:2}
Let $(\lambda_i)$ be a  solution of a Garnier system governing the isomonodromic deformation of a trace free logarithmic connection $\nabla$ on $\mathcal O_{\Pu}^{\oplus 2}$  with no apparent pole.
The following are equivalent.
\begin{enumerate}
\item \label{appli1intro}The multivalued functions $\lambda_i$  are algebraic functions.
\item \label{appli2intro}The functions $\lambda_i$ have finitely many branches.
\item \label{appli3intro} The conjugacy class $[\rho]$ of the monodromy representation of $\nabla$ has finite orbit under $\mathrm{MCG}_n\Pu$.
\end{enumerate}
\end{THM}

For $n\geq 4$ punctures, there exist irreducible rank $m$ representations that yield finite orbits, so that the above statements are not vacuous; see section~$\ref{non vacuity}$.

~\\
We now want to comment on the possible continuations of this work, after we will explain the structure of the paper.

~\\

The most interesting part of Theorem \ref{maintheorem} is the implication $(\ref{unthmintro})\Rightarrow (\ref{deuxthmintro})$, it allows to construct flat connections on vector bundles over projective ruled varieties from finite orbits.  We already have an example that testifies that, in general, the obtained connection is not too simple, see \cite{hilbmod}.
The determination of conjugacy classes of representations that have finite orbits for a fixed rank $m$ and number of punctures $n$ seems to be a good strategy to exhibit interesting flat connections.
For $(m,n)=(2,4)$ this study was initiated by Dubrovin-Mazzocco \cite{MR1767271}, continued by Iwasaki \cite{MR1930217}, Boalch \cite{MR2254812}, Cantat-Loray \cite{MR2649343} and completed by  Lisovyy-Tykhyy \cite{MR3253555}. 

For $(m,n)=(2,4)$ this study was initiated by Dubrovin-Mazzocco \cite{MR1767271}, systematized by Cantat-Loray \cite{MR2649343} and completed by  Lisovyy-Tykhyy \cite{MR3253555}. Thanks to the efforts of many people, such as Hitchin \cite{MR1351506}, Dubrovin-Mazzocco \cite{MR1767271}, Doran \cite{MR1909248}, Andreev-Kitaev \cite{MR1911252}, Boalch \cite{MR2681697}, the algebraic solutions of the corresponding Garnier system (Painlev\'e VI) are essentially all known.

For $m=2$, $n>4$, restricting to the case of Zariski dense representations in $\mathrm{SL}_2(\C)$, Diarra \cite{MR3077637} has determined all the finite orbits that can be obtained \textit{via} the  pullback technique used in \cite{MR1909248} and \cite{MR1911252}. They are finite and quite few. One may also find algebraic Garnier solutions directly from algebraic Painlev\'e VI solutions, as illustrated in \cite{GirandGarnier}. 
It might be interesting to determine the other finite orbits, say for $n=5$. Then, a challenging program would be to exhibit the corresponding algebraic Garnier solutions. 

We should also point out that we have an important structure theorem for Zariski dense rank two representations of quasi-projective fundamental groups given by Corlette and Simpson in \cite{MR2457528}. For rank two reducible representations, this theorem is well complemented by  \cite[Theorem 5.1]{Bartolo:arXiv1005.4761}. Both results are fruitfully used in the theory of foliations, see \cite{MR3294560} and \cite{loray2014representations}.

To the authors knowledge, in the general quasi-projective case, there is no known similar result for higher rank. In this regard, investigation of the finite orbits occurring in Theorem~\ref{maintheorem} with $m>2$ could be an interesting testing ground for possible generalization of \cite{MR2457528}. In section $\ref{Fields}$, we give some information on the fields of definitions of the projective representations that would occur in these examples.

Other interesting questions concern the generalization of Theorem \ref{maintheorem} for isomonodromic deformations of logarithmic connections over other Riemann surfaces.
Also, in the side of representations, an important part of our construction is not restricted to representations with value in $\mathrm{GL}_m(\C)$ and we may ask if we could exploit it in a context of isomonodromic deformations of connections on more general principal bundles.

~\\

The paper is organized as follows.
The proof of Theorem \ref{maintheorem} is based on  a form of Riemann-Hilbert correspondence. In section $\ref{topconstruction}$, from a finite orbit, we build a representation of the fundamental group of the complement of a hypersurface in a projective ruled variety. In section $\ref{secRH}$, we prove the appropriate Riemann-Hilbert correspondence. We combine these results in the beginning of section $\ref{isomdef}$ and obtain Theorem \ref{maintheorem}. In section $\ref{secGarnier}$,
we  introduce the relation between isomonodromic deformations and Garnier systems and conclude by the proof of Theorem \ref{COR:2}. 
The paper is concluded by section \ref{Fields}, a discussion on fields of definitions of projective representations of semi-direct products, whose interest has been evoked above.
A technical point on Garnier systems is treated in the Appendix.
~\\

{\bf Acknowledgements.}
The first topological ideas about this work were found whilst the author worked at IMPA, there we wish to thank  A. Lins Neto,  J.V. Pereira and  P. Sad  for useful discussions.
The main part of this work was accomplished in  the Mathematics department of the University of Pisa. We acknowledge M. Salvetti and  F. Callegaro for discussions on hyperplane arrangements and lifting issues. The author was hosted in Pisa by M. Abate and J. Raissy with a postdoc fellowship in the framework of Italian FIRB project  \textit{Geometria Differenziale e Teoria Geometrica delle Funzioni}.
We are grateful  for this hospitality and  this funding, as for an additional contribution of J. Raissy through some helpful proof reading.  We thank both referees for their relevant suggestions and remarks. Last but not least, the author would like to thank F. Loray for introduction to the topic of isomonodromic deformations and for numerous enlightening conversations.

\section{A topological construction from braid group orbits } \label{topconstruction}
\subsection{Braid group and Mapping class group}\label{BraidMCG}
 We recall here some known results about mapping class groups, braid groups and their relations. The main reference is \cite{MR0243519}.

For any connected real manifold $M$, let $\mathrm{Homeo}(M)$ denote the group of self homeomorphisms of $M$. Let $b\in M$ be a base point. For $h\in\mathrm{Homeo}(M)$, we have an isomorphism $h_*: \pi_1(M,b)\rightarrow \pi_1(M,h(b))$, for any choice of path $\gamma$ in $M$ from $b$ to $h(b)$, we have an isomorphism $\gamma_*:\pi_1(M,h(b))\simeq \pi_1(M,b) $, and $\gamma_*\circ h_*$ is an element of $\mathrm{Aut}\left(\pi_1(M,b)\right)$, the group of automorphisms of $\pi_1(M,b)$. If we change the path $\gamma$, the class of $\gamma_*\circ h_*$ in the group $\mathrm{Out}\left(\pi_1(M,b)\right)$ of outer automorphisms does not change.
We get a morphism  $$\M(M)\rightarrow \mathrm{Out}\left(\pi_1(M,b)\right).$$
Here $\M(M)$ is the mapping class group of $M$, that is the group of isotopy classes in $\mathrm{Homeo}(M)$. In our work we are especially interested in the $n$ times punctured Riemann sphere $M=\Pu(\C) \setminus \{x_1,\ldots,x_n\}$ and the $n!$-index subgroup of $\M(M)$ which fixes each of these punctures. We denote this subgroup by $\PM_n(\Pu)$, the letter $\mathrm{P}$ standing for \emph{pure} mapping class group, in opposition to the \emph{full} mapping class group  $$\M_n(\Pu):=\M(\Pu(\C)\setminus \{x_1,\ldots,x_n\}).$$ 

In this special setting, for $n>2$, the map $$\PM_n(\Pu)\rightarrow \mathrm{Out}\left(\pi_1(\Pu(\C)\setminus \{x_1,\ldots,x_n\})\right)$$ can be recovered from an action of the pure braid group on $n$ strands on the sphere. This will be useful for our considerations and
we explain it in the next paragraphs.

Let $F_{0,n} \Pu:=\{(y_1,\ldots,y_n)\in (\Pu)^n \vert  i\neq j \Rightarrow y_i\neq y_j\}.$
The pure braid group on $n$ strands on the sphere is $\pi_1(F_{0,n}\Pu)$.

For $n>2$, consider $$F_{3,n-3}\Pu:=\{(t_1,\ldots,t_{n-3},0,1,\infty)\in F_{0,n}\Pu\}.$$

Normalizing configurations, we get a homeomorphism $$\begin{array}{ccccc} \mathrm{PSL}_2(\C) &\times& F_{3,n-3}\Pu& \simeq& F_{0,n}\Pu\\
																							\left[~A \right.&,& (t_1,\ldots, t_n\left.)~\right]& \rightarrow &(A\cdot t_1,\ldots,A\cdot t_n)\\ \end{array}$$

We have a natural projection $$p_{n+1} : F_{0,n+1}\Pu \rightarrow F_{0,n}\Pu, (y_0,\ldots,y_{n})\mapsto (y_1,\ldots,y_{n}).$$ 
This projection $p_{n+1}$ has a continuous section.
Namely, we define such a section by $$\sigma_n(A\cdot(t_1,\cdots,t_{n-3},0,1,\infty))=A\cdot(\tau(t),t_1,\cdots,t_{n-3},0,1,\infty),$$ $$\mbox{with }  \tau(t)=2+\sum_{i=1}^{n-3} \Re(t_i)^2,$$
where $\Re(t_i)$ stands for the \emph{real part} of the complex number $t_i$.

We denote $s_n(y)$ the first coordinate of $\sigma_n(y)$. In this way we have
\[\sigma_n(y)=(s_n(y),y_1,\ldots,y_n),~~~~y\in F_{0,n}\Pu.\]

Fix a base point $x=(x_1,\ldots,x_n)\in F_{0,n}\Pu$ and let $x_0=s_n(x).$
The group $\pi_1(\Pu\setminus \{x_1,\ldots,x_n\},x_0)$ will be ubiquitous in this paper. For brevity of notation we define, once and for all,
$$\gf:=\pi_1\left(\Pu\setminus \{x_1,\ldots,x_n\},x_0\right).$$

The section $\sigma_n$ induces a morphism $${\sigma_{n}}_* : \pi_1(F_{0,n}\Pu,x)\rightarrow\pi_1(F_{0,n+1}\Pu,\sigma_n(x)).$$

For $\M_{n+1}(\Pu)=\M(\Pu \setminus \{x_0,\ldots,x_{n}\})$, \cite{MR0243519} gives two other maps
$$\left \lbrace \begin{array}{l}
 d_*: \pi_1(F_{0,n+1}\Pu)\rightarrow \mathrm{P}\M_{n+1}(\Pu)\mbox{, [Th.$1$ p.$216$] and }\vspace{0.1cm}\\
 \mu_{n+1} :\mathrm{PMCG}_{n+1}(\Pu)\rightarrow \mathrm{Aut}(\gf)\mbox{, [sec.$3$]}.
\end{array} \right.$$

 The action of $\mathrm{P}\M_{n+1}$ through  $\mu_{n+1}$ is simply by composition on the level of representative loops: $(g\cdot\alpha)(t):=g(\alpha(t))$.
The map $d_*$ can be described as follows. For $\alpha\in \pi_1(F_{0,n}\Pu,x)$, let $(\alpha_i(t))_i : [0,1] \rightarrow F_{0,n}\Pu\subset \left(\Pu\right)^n$ be a representative loop. There exists a continuous family $h_t\in \mathrm{Homeo}(\Pu)$ such that $h_t(x_i)=\alpha_i(t)$.  The class $[h_1]$ of $h_1$ in $\M_n(\Pu)$ depends only on $\alpha$ and we have $d_*\alpha=[h_1] \in \M_n(\Pu)$. We see that, with the group law induced by composition of maps on  $\M_n$, $d_*$ is an \emph{antimorphism}. However, $\mu_{n+1}$ is a genuine morphism.

It is readily checked that, for $n>2$, we have the following commutative diagram.
 $$\xymatrix{
 \pi_1(F_{0,n}\Pu,x) \ar[d]^{d_*}\ar[r]^{\sigma_{n*}} & \pi_1(F_{0,n+1}\Pu,\sigma_n({x}))\ar[d]^{d_*}\\
   \mathrm{P}\M_{n}(\Pu)  \ar[d]& \mathrm{P}\M_{n+1}(\Pu)\ar[l]   \ar[d]^{\mu_{n+1}}\\
 \mathrm{Out}\left(\gf\right)&\ar[l] \mathrm{Aut}\left(\gf\right)    
  }$$
  This is the announced relation between braid group and mapping class group actions.
  
At this point, we should notice a few more facts. The kernel of $d_*:\pi_1(F_{0,n}\Pu,x) \rightarrow \PM_{n}(\Pu)$ is contained in the center of $\pi_1(F_{0,n}\Pu)$,  \cite[Cor.$1.2$ p.$218$]{MR0243519}.
The center of $\pi_1(F_{0,n}\Pu)$ is  is isomorphic to $\Z/2\Z\simeq\pi_1(\mathrm{PSL}_2(\C))$, see \cite[Prop.$7$ p.$311$]{MR2092062}.
One checks directly that the image of $\pi_1(\mathrm{PSL}_2(\C))$ \textit{via} $\mathrm{PSL}_2(\C) \times F_{3,n-3}\Pu \simeq F_{0,n}\Pu$ is contained in this kernel.

 Also, we have an exact sequence \cite[Th.$1$ p.$216$]{MR0243519}
$$\pi_1(F_{3,n-3}\Pu) \stackrel{d_*}{\rightarrow} \PM_{n}(\Pu)\rightarrow \PM_3(\Pu),$$ which asserts surjectivity of $d_*:\pi_1(F_{0,n}\Pu) \rightarrow \PM_{n}(\Pu)$, because $\PM_3(\Pu)$ is the trivial group.
Finally, $d_*$ induces an anti-isomorphism $$\pi_1(F_{3,n-3}\Pu)\simeq  \PM_n(\Pu),~~~n>2.$$

The map $\psi :\pi_1(F_{0,n}\Pu,x)\rightarrow \mathrm{Aut}\left(\gf\right)$ induced by the above diagram will be of major importance in the next section.
Occasionally, we will denote $\alpha^\beta:=\psi(\beta)\cdot \alpha$; for $\alpha\in \gf$, $\beta \in \pi_1(F_{0,n}\Pu,x)$.
\subsection{Fundamental groups of locally trivial fiber bundles}\label{groupeloctriv}
We discuss fundamental groups of certain locally trivial fiber bundles. Our treatment is inspired by \cite{shimada}, to which we refer for introduction to braid monodromy.
\begin{theorem}[Zariski-Van Kampen]\label{ZarVan}
Let $p: E\rightarrow B$ be a topologically locally trivial fiber bundle. Suppose $E$ is path-connected and we have a continuous section $s: B \rightarrow E$, $p\circ s=id_{B}$. denote $E_b:=p^{-1}(b)$, for a base point $b\in B$.
  Then, the fundamental group of $E$ is a semidirect product:
$$\pi_1\left(E,s(b)\right)\simeq\pi_1(E_b,s(b))\rtimes_{\phi}\pi_1(B,b);$$
 
More precisely, the right-hand side is the product $\pi_1(E_b,s(b))\times\pi_1(B,b)$ endowed with multiplication  law 
$$(n_1,h_1)\star(n_2,h_2)=\left( n_1(\phi(h_1)\cdot n_2),h_1 h_2 \right),$$ 
for a morphism $\phi :\pi_1(B,b)\rightarrow \mathrm{Aut}(\pi_1(F_b,s(b)))$  characterized by the following property:

 for any $(\alpha,\beta) \in \pi_1(E_b,s(b))\times\pi_1(B,b)$, $$\phi(\alpha) \cdot \beta=s_*(\alpha)\beta s_*(\alpha)^{-1}\mbox{ in } \pi_1(E,s(b)).$$
\end{theorem}
\begin{proofbis}
We have a long homotopy exact sequence.
$$\cdots\rightarrow\pi_2(E,s(b))\stackrel{p_*}{\rightarrow}\pi_2(B,b) \rightarrow \pi_1(E_b,s(b))\stackrel{i_*}{\rightarrow}\pi_1(E,s(b))\stackrel{p_*}{\rightarrow}\pi_1(B,b)\rightarrow \{1\}$$

The section $s$ induces a section $s_*$ of $\pi_2(E,s(b))\stackrel{p_*}{\rightarrow}\pi_2(B,b)$, by exactness we deduce the short exact sequence
$$\{1\}\rightarrow \pi_1(E_b,s(b))\stackrel{i_*}{\rightarrow}\pi_1(E,s(b))\stackrel{p_*}{\rightarrow}\pi_1(B,b)\rightarrow \{1\}$$ which splits by $s_*$, yielding the announced semidirect product structure.
\qed \end{proofbis}
We are specially interested with the following fiber bundles.
\begin{definition} \label{tautbundle}
For $n>2$, let $$\mathcal{T}_n:=\{(y,z)\in F_{0,n}\Pu\times \Pu \vert y_i\neq z, i=1,\ldots,n\}.$$ The natural map $\mathcal{T}_n\rightarrow F_{0,n}\Pu$ is called  \emph{the tautological fiber bundle on~$F_{0,n}\Pu$}.
 We endow this bundle with the section $s_n$ introduced in section $\ref{BraidMCG}$.

For any subgroup $H$ of $\pi_1(F_{0,n}\Pu,x)$, choose a covering $f_H: (B_H,b_H)\rightarrow (F_{0,n}\Pu,x)$ with ${f_H}_*\pi_1(B_H,b_H)=H$.  We define $(E_H,s_H)$ as the pullback of $(\mathcal{T}_n,s_n)$ by $f_H$.

\end{definition}

\begin{proposition} \label{semidirectsphereU}
Let $f : (B,b)\rightarrow (F_{0,n}\Pu,x)$ be a continuous map, with $B$ path connected.
By pullback, we have a commutative diagram 
$$\xymatrix@!{
 f^*\mathcal{T}_n\ar[d] \ar[r] & \mathcal{T}_{n} \ar[d]\\
B\ar[r]^f \ar@/^1.3pc/[u]^{s_f}&F_{0,n}\Pu \ar@/_1.3pc/[u]_{s_n}
}$$
We have a semidirect product decomposition
$$\pi_1\left({f^*\mathcal{T}_n},s_f(b)\right)\simeq\gf\rtimes_{\varphi}\pi_1(B,b)$$ with structure map
$$\varphi=\iota\circ \psi \circ f_*$$ where 
\begin{itemize}
 \item ${f}_* : \pi_1(B,b)\rightarrow \pi_1(F_{0,n}\Pu,x)$ is the map induced by $f$,
 \item  $\psi :\pi_1(F_{0,n}\Pu,x)\rightarrow \mathrm{Aut}\left(\gf\right)$ is the antimorphism introduced in section $\ref{BraidMCG}$.
 \item  $\iota$ is the involution $g\mapsto g^{-1}$.
  \end{itemize}
\end{proposition}
\begin{proofbis}  In the proof we denote $E:=f^*\mathcal{T}_n$, $s:=s_f$ and $E_b$ the fiber of $b$ in $E$.

First, $\mathcal{T}_n$ and $E$ are topologically locally trivial fiber bundles.
The local triviality of $\mathcal{T}_n$ is an exercise which reduces to finding a continuous family $(h_{\tau})_{\tau \in \mathbb{D}}$ of homemorphisms of the closed disk $\overline{\mathbb{D}}$ with $h_{\tau}(\tau)=0$. 
 In the same way, we can prove that $E\setminus s(B)\rightarrow B$ is also locally trivial.

  By Theorem $\ref{ZarVan}$, we have a semidirect product structure $$\pi_1\left(E,s(b)\right)\simeq\pi_1(E_b,s(b))\rtimes_{\phi}\pi_1(B,b);$$
  let us characterize the structure map $\phi$.
 Take loops $\alpha : [0,1]\rightarrow (B,b), \beta : [0,1]\rightarrow \left(E_b,s(b)\right)$.
  Through local trivializations of $E\setminus s(B)$ over $\alpha$ we can define a continuous family $(h_t)_{t \in [0,1]}$ of homemorphisms $h_t: E_{b}\setminus s(b) \stackrel{\simeq}{\longrightarrow} E_{\alpha(t)}\setminus s(B)$, with $h_0=id_{E_b}$.
  In terms of the projections $(y,z)$ to the factors of the product $B \times \Pu$, $h_t$ takes the form $h_t: (b,z)\mapsto (\alpha(t),g_t(z))$ for a homeomorphism $$g_t: \Pu\setminus \{x_0,\ldots,x_{n}\} \stackrel{\simeq}{\longrightarrow}  \Pu \setminus \cup_{i=1}^{n+1} \{pr_i\circ \sigma_n\circ f\circ\alpha(t)\},$$
  where $pr_i : F_{0,n+1}\Pu \hookrightarrow\left(\Pu\right)^{n+1} \rightarrow \Pu$ is the projection to the $i$-th factor.
  Hence, the resulting element $[g_1]\in \mathrm{P}\M_{n+1}(\Pu)$ is $d_*{\sigma_n}_*f_*\alpha$ and the corresponding element of $\mathrm{Aut}(\pi_1(\Pu\setminus\{x_1,\ldots,x_n\},x_{0})$ is $\psi ( f_*\alpha)$.
 
 The continuous family  $\beta_t(s):=h_t(\beta(s)), t\in[0,1]$ is such that $\beta_t$ is a loop in $E_{\alpha(t)}$ with base point $s(\alpha(t))$ and $\beta_0=\beta$. Thanks to this, we easily see that the loop $\beta_1$ in $E_b$ represents
$s_*(\alpha)^{-1}\beta s_*(\alpha)$ in $\pi_1(E,s(b))$, thus it represents $\phi(\alpha)^{-1}\cdot \beta$ in $\pi_1(E_b,b)$. We just proved $\phi(\alpha)^{-1}\cdot \beta=h_{1*}\beta$.

We now identify $(E_b,b)$ with $(\Pu\setminus \{x_1,\ldots,x_n\},x_{0})$ \textit{via} the projection~$z$. This induces  an identification $z_*: \pi_1(E_b,b)\stackrel{\simeq}{\longrightarrow}\pi_1(\Pu\setminus \{x_1,\ldots,x_n\},x_{0})=\gf$. 
Thanks to this, we also have a semidirect product structure $$\pi_1\left(E,s(b)\right)\simeq\gf\rtimes_{\varphi}\pi_1(B,b),$$ and the identity $\phi(\alpha)^{-1}\cdot \beta=h_{1*}\beta$ reads $\varphi(\alpha)^{-1}\cdot z_*\beta=g_{1*} \cdot z_*\beta$, which implies, by our previous description of ${g_1}_*$,
$\varphi(\alpha)^{-1}=\psi (f_*\alpha)$. We have the announced group isomorphism.
\qed \end{proofbis}
\begin{corollary}\label{semidirectsphere}
Let $\rev_H :(B_H,b_H)\rightarrow (F_{0,n}\Pu,x)$ and $E_H$ be as in Definition $\ref{tautbundle}$.

We have a semidirect product decomposition 
$$\pi_1\left(E_H,s_H(b)\right)\simeq\gf\rtimes_{\zeta}H $$ with structure map
$\zeta=\iota\circ \psi_{\vert H}$, where $\iota$ is the involution $g\mapsto g^{-1}$.
\begin{proofbis} We identify $\pi_1(B_H,b_H)$ with $H$ \textit{via} ${f_H}_*$.
\qed \end{proofbis}
\end{corollary}
We want to use  Corollary $\ref{semidirectsphere}$ to give informations on linear representations of $\pi_1(E_H)$. Yet, the first informations we get are about their projectivizations.
\subsection{Projective representations}
With the notation of section \ref{groupeloctriv}, we will study a group morphism $\rho: \pi_1(E_H,s_H(b_H))\rightarrow G$  from the properties of its restriction $\rho_b: \gf\rightarrow G$ to the first factor in the semidirect product structure of Corollary $\ref{semidirectsphere}$.

We denote $\mathrm{P}$ the projection map $\mathrm{P}:\mathrm{GL}_m(k)\rightarrow \mathrm{GL}_m(k)/k^*$.
Take a representation $\rho:\pi_1(E_H,s_H(b_H))\rightarrow \mathrm{GL}_m(k)$, for $k$ an algebraically closed field. Suppose $\rho_b$ is irreducible, by Schur's lemma the image of
$$\mathrm{P}\rho_b: \gf\rightarrow\mathrm{P}\mathrm{GL}_m(k)$$ has trivial centralizer in $\mathrm{P}\mathrm{GL}_m(k)$ and the following applies, showing the projectivization of $\rho$ is characterized by $\rho_b$ in a quite effective manner.
\begin{lemma}\label{lemsemidirect}
Consider the following split exact sequence of groups.
    $$ \xymatrix{
    \{1\} \ar[r]&  N \ar[r]^{i}& \Gamma \ar[r]& Q \ar@/^-1pc/[l]_{s}  \ar[r]&\{1\}
    }$$
Let $G$ be a group. Suppose we have a morphism $\rho:i(N)\rightarrow G$, such that $\mathrm{Im}(\rho)$ has trivial centralizer in $G$. Then $\rho$ extends to at most one morphism $\hat{\rho}: \Gamma \rightarrow G$.

 Moreover, the existence of such an extension is equivalent to the condition
 \begin{equation}\label{stabilisation} \exists (g_{\beta}) \in G^{Q}\mbox{ such that }\rho(\beta^{-1} \alpha \beta)=g_{\beta}^{-1}\rho(\alpha)g_{\beta}~;~\forall (\alpha,\beta)\in i(N)\times s(Q).
 \end{equation}
 In addition, if such a family $(g_\beta)$ exists, it is unique.
\begin{proofbis}
We must have, for any $(\alpha,\beta) \in i(N) \times s(Q)$, 
$$\rho(\beta^{-1} \alpha \beta)=\hat{\rho}(\beta)^{-1}\rho(\alpha)\hat{\rho}(\beta).$$

 For $i=1,2$, let $g_{\beta,i}$ statisfy $$\rho(\beta^{-1} \alpha \beta)=g_{\beta,i}^{-1}\rho(\alpha)g_{\beta,i}~;~~~\forall \alpha\in i(N).$$
 Then $g_{\beta,1}g_{\beta,2}^{-1}$ centralizes $\mathrm{Im(\rho)}$ and must be trivial.
 
 Thus if $(g_\beta)$ exists, it is unique and $g_{\beta_1\beta_2}=g_{\beta_1}g_{\beta_2}$.
 Then one checks that the assignment $\alpha\beta\mapsto \rho(\alpha)g_\beta$ defines the required extension $\hat{\rho}$.
\qed \end{proofbis}

\end{lemma}
In the next section, we will see how condition $(\ref{stabilisation})$ appears naturally in our study.
\subsection{Braid group action on the set of representations}
In section $\ref{BraidMCG}$ we have described an antimorphism 
$$\psi: \pi_1(F_{0,n}\Pu,(x_1,\ldots,x_n))\rightarrow \mathrm{Aut}\left(\gf\right).$$

It induces an action $(\beta,\rho)\mapsto \beta\cdot\rho$ of $\pi_1(F_{0,n}\Pu,(x_1,\ldots,x_n))$ on the set of morphisms $\mathrm{Hom}\left(\gf,G\right)$, for any group $G$. Explicitly: $(\beta\cdot \rho)(\alpha):=\rho(\alpha^\beta)$. This action commutes with the action of $G$ on $\mathrm{Hom}\left(\gf,G\right)$ by inner automorphisms and we have an induced action $(\beta,[\rho])\mapsto [\beta\cdot\rho]$; where $[\rho]$  denotes the class in $\mathrm{Hom}(\gf,G)/G$ of the element $\rho \in \mathrm{Hom}(\gf,G)$.

In this regard, in the case where $\Gamma=\pi_1(f^*\mathcal{T}_n, s_f(b))$ as in Proposition $\ref{semidirectsphereU}$, condition $(\ref{stabilisation})$  of Lemma $\ref{lemsemidirect}$ means that $f_*\pi_1(B,b)<\pi_1(F_{0,n}\Pu,x)$ is contained in the stabilizer $\mathrm{Stab}[\rho]$ of $[\rho]\in\mathrm{Hom}(\gf,G)/G$, independently of any assumption on $G$ or $\rho$.
We have straightforwardly three corollaries of this observation and Lemma $\ref{lemsemidirect}$.

\begin{corollary}\label{extproj}
Consider a morphism $\rho: \gf\rightarrow \mathrm{PGL}_m(k)$, for $k$ a field.

Suppose $\mathrm{Im}(\rho)$ has trivial centralizer in $\mathrm{PGL}_m(k)$ and $H<\pi_1(F_{0,n}\Pu,x)$ stabilizes $[\rho]\in\mathrm{Hom}\left(\gf,\mathrm{PGL}_m(k)\right)/\mathrm{PGL}_m(k)$.

Then $\rho$ is the restriction $\hat{\rho}_{b_H}$ of a unique representation $$ \hat{\rho}: \pi_1(E_{H},s_H(b_H))\rightarrow \mathrm{PGL}_m(k). $$\qed
\end{corollary}

\begin{corollary}\label{corlin}
Let $k$ be an algebraically closed field and
 $\rho: \gf \rightarrow \mathrm{GL}_m(k)$ be a group morphism.

Suppose $\rho$ is irreducible and $H<\pi_1(F_{0,n}\Pu,x)$ stabilizes $$[\rho]\in\mathrm{Hom}\left(\gf,\mathrm{GL}_m(k)\right)/\mathrm{GL}_m(k).$$

Then $\mathrm{P}\rho$, the projectivization of $\rho$,  is the restriction ${\widehat{\mathrm{P}\rho}}_b$  of a unique projective representation $$\widehat{\mathrm{P}\rho}: \pi_1(E_{H},s_H(b_H))\rightarrow \mathrm{PGL}_m(k).$$
\qed
\end{corollary}
\begin{corollary}\label{finiteindex}
Let $f: (B,b)\rightarrow (F_{0,n}\Pu,x)$ be a continuous map, with $B$ path connected.
Consider a group morphism $$\rho: \pi_1(f^*\mathcal{T}_n,s_f(b))\simeq\gf\rtimes_{\varphi}\pi_1(B,b) \rightarrow G.$$
If   $f_*: \pi_1(B,b)\rightarrow\pi_1(F_{0,n}\Pu,x)$ has finite index image and $\rho_b: \gf\rightarrow G$ is the restriction of $\rho$ to $\gf$, then the element $$[\rho_b]\in \mathrm{Hom}(\gf,G)/G$$ has a finite orbit under the action of $\pi_1(F_{0,n}\Pu,x)$.
\qed
\end{corollary} 
\subsection{Lifting issues}
A natural question after stating Corollary $\ref{corlin}$ is the following.
\begin{question} \label{queslift} Does the extension $\widehat{\mathrm{P}\rho}$ constructed \textit{via} Corollary $\ref{corlin}$ lift to a linear extension of $\rho$ ?
\end{question}

Recall one says a representation  $\mu :\Gamma \rightarrow \mathrm{PGL}_m(k)$ lifts to $\mathrm{GL}_m(k)$   if we can find a map $\tilde{\mu}: \Gamma \rightarrow \mathrm{GL}_m(k)$ such that $\mu=\mathrm{P}\tilde{\mu}$.
If we can find a lift $\tilde{\mu}$  with image in $\mathrm{SL}_m(k)$, we say that $\mu$ lifts to $\mathrm{SL}_m(k)$.
\begin{lemma}\label{lemliftabstr} 
Let $r: N\rtimes Q \rightarrow \mathrm{PGL}_m(k)$ be a group morphism, for a given field $k$.

Suppose the restrictions of $r$ to both factors, $r_N$ and $r_Q$,  lift to $\mathrm{GL}_m(k)$ and $\mathrm{Im}(r_N)$ has trivial centralizer in $\mathrm{PGL}_m(k)$.
Suppose a lift $\tilde{r}_N$ of $r_N$ satisfies \begin{equation}\label{stabilisationbis} \exists (g_{\beta}) \in \left(\mathrm{GL}_m(k)\right)^{Q} \vert~~\forall (\alpha,\beta)\in N\times Q,~ \tilde{r}_N(\beta^{-1} \alpha \beta)=g_{\beta}^{-1}\tilde{r}_N(\alpha)g_{\beta}.
 \end{equation}
 Then $r$ lifts to $\mathrm{GL}_m(k)$. If $r_N$ and $r_Q$ lift to $\mathrm{SL}_m(k)$, then so does $r$.
\end{lemma}
\begin{proofbis}
Let $\tilde{r}_N,\tilde{r}_Q$ be lifts for $r_N$ and $r_Q$, with values in $\mathrm{SL}_m(k)$ if possible.
The relations of condition $(\ref{stabilisationbis})$ projectivize to $r_N(\beta^{-1} \alpha \beta)=\mathrm{P}g_{\beta}^{-1}r_N(\alpha)\mathrm{P}g_{\beta}.$

If we replace $\mathrm{P}g_{\beta}$ by $r_Q(\beta)$, these equations still hold, because $r$ is a morphism. By Lemma $\ref{lemsemidirect}$, we get $\mathrm{P}g_{\beta}=r_Q(\beta)$ and there exists a family $(\lambda_\beta)$ of non zero scalars such that
$g_{\beta}=\lambda_{\beta}\tilde{r}_Q(\beta)$, in particular $$\tilde{r}_N(\beta^{-1} \alpha \beta)=\tilde{r}_Q(\beta)^{-1}\tilde{r}_N(\alpha)\tilde{r}_Q(\beta).$$ and the assignment $\alpha\beta\mapsto \tilde{r}_N(\alpha)\tilde{r}_Q(\beta)$ defines the desired lift of  
$r$.
\qed \end{proofbis}
\begin{corollary}\label{corlift} 
Let $\rho\in \mathrm{Hom}\left(\gf,\mathrm{GL}_m(k)\right)$, $k$ an algebraically closed field.
Suppose $\rho$ is irreducible and $H<\pi_1(F_{0,n}\Pu,x)$ stabilizes $$[\rho]\in\mathrm{Hom}\left(\gf,\mathrm{GL}_m(k)\right)/\mathrm{GL}_m(k).$$
Consider  $\widehat{\mathrm{P}\rho} : \pi_1(E_H,s_H(b_H)) \rightarrow \mathrm{PGL}_m(k)$, the  extension of $\mathrm{P}\rho$ obtained in Corollary $\ref{corlin}$.
Let $g : (B,b)\rightarrow (B_H,b_H)$ be a continuous map, with $B$ path connected and let $f=f_H\circ g$.
Let $$\widehat{P\rho}_{\vert B} : \pi_1(f^*\mathcal{T}_n, s_{f}(b))\simeq \gf\rtimes_{\varphi}\pi_1(B,b)\rightarrow \mathrm{PGL}_m(k)$$ be the pullback of $\widehat{P\rho}$ by the natural map $f^*\mathcal{T}_n\rightarrow E_H$.

Suppose,   the restriction $r_{\pi_1(B,b)}$ of $\widehat{\mathrm{P}\rho}_{\vert B}$ to the second factor $\pi_1(B,b)$ lifts to $\mathrm{GL}_m(k)$.

 Then $\widehat{\mathrm{P}\rho}_{\vert B}$ lifts to $\mathrm{GL}_m(k)$ as an extension of $\rho$. If $\rho$ takes values in $\mathrm{SL}_m(k)$, and $r_{\pi_1(B,b)}$ lifts to $\mathrm{SL}_m(k)$, then $\widehat{\mathrm{P}\rho}_{\vert B}$ lifts to $\mathrm{SL}_m(k)$ as an extension of $\rho$.
\end{corollary}
\begin{proofbis}
Take $r=\widehat{\mathrm{P}\rho}_{\vert B}$, $\tilde{r}_N=\rho$ in Lemma $\ref{lemliftabstr}$. The condition $(\ref{stabilisationbis})$ is ensured by $H$ stabilizing $[\rho]$.
 \qed \end{proofbis}
\begin{corollary}
In the case $n=4$ of Corollary $\ref{corlin}$, $\rho$ is in fact the restriction $\hat{\rho}_b$ of a linear representation  $\hat{\rho}:\pi_1(E_{H},s_H(b_H))\rightarrow \mathrm{GL}_m(k).$
If $\rho$ takes values in $\mathrm{SL}_m(k)$, so can be chosen $\hat{\rho}$.
\begin{proofbis} We will use Corollary $\ref{corlift}$ for $B=B_H$, $g=id_B$. We want to lift the restriction $r_{\pi_1(B,b)}=r_H$ to $\mathrm{SL}_m(k)$.
We have the product structure $\pi_1(F_{0,4}\Pu)\simeq\pi_1(\mathrm{PSL}_2(\C))\times\pi_1(F_{3,1}\Pu)=\mathbb{Z}/2\mathbb{Z}\times F_2$, where $F_{\ell}$ is the free group of rank $\ell$. We know the action of $\pi_1(\mathrm{PSL}_2(\C))$ on $\mathrm{Hom}(\pi_1(\Pu\setminus \{x_1,\ldots,x_n\}),\mathrm{GL}_m(k))$
is trivial, because the section $s_n$ is $\mathrm{PSL}_2(\C)$-equivariant.

 Thus if $H=\mathrm{Stab}[\rho]$, then $H=\mathbb{Z}/2\mathbb{Z}\times F_{\ell}$ for some $\ell$ and the representation $r_H$ is trivial on the factor $\mathbb{Z}/2\mathbb{Z}$. For this reason, it is enough to lift the restriction of $r_H$
 to the factor $F_{\ell}$, which is easy because $F_{\ell}$ is free. Then if we take for $H$ a subgroup of $\mathrm{Stab}[\rho]$, the restriction of the lift we just constructed for $\mathrm{Stab}[\rho]$ gives the required lift of $r_H$.
\qed \end{proofbis}
\end{corollary}
\begin{remark}
The linear extension of $\rho$ obtained in Corollary $\ref{corlift}$ is almost never unique, because we can tensor any lift by a scalar representation of $\pi_1(B,b)$ to obtain a new lift.
\end{remark}

Corollary $\ref{corlift}$ shows Question $\ref{queslift}$ reduces to determining if a certain projective representation of a subgroup $H$ of the pure braid group 
lifts to a linear representation.

One can try to solve this issue using cohomological arguments. This takes the following form.
First, projective and linear representations of $H$ correspond to flat bundles on $B_H$, this means elements of the cohomology sets $\mathrm{H}^1(B_H,\mathrm{PGL}_m(k))$ and $\mathrm{H}^1(B_H,\mathrm{GL}_m(k))$ respectively. 

As it is the application we have in mind, we restrict to the case $k=\C$. We have a commutative diagram of locally constant sheaves, the lines of which are exact sequences.
 $$\xymatrix{
0\ar[r]&\mathbb{Z}/m\mathbb{Z} \ar[r] \ar[d]^{i}&\mathrm{SL}_m(\C) \ar[r] \ar[d]^{i}&\mathrm{PGL}_m(\C) \ar[r]\ar@{=}[d]&0\\
 0\ar[r]&\C^*\ar[r]&\mathrm{GL}_m(\C)\ar[r]&\mathrm{PGL}_m(\C)\ar[r]&0\\
  }$$
This yields a commutative diagram of cohomology sets, with exact sequences as lines.
$$\xymatrix{  \mathrm{H}^1(B_H,\mathrm{SL}_m(\C))\ar[d] \ar[r] & \mathrm{H}^1(B_H,\mathrm{PGL}_m(\C)) \ar[r]^{d_1} \ar@{=}[d]&\mathrm{H}^2(B_H,\mathbb{Z}/m\mathbb{Z} )\ar[d]^{i_*}\\
\mathrm{H}^1(B_H,\mathrm{GL}_m(\C))\ar[r]& \mathrm{H}^1(B_H,\mathrm{PGL}_m(\C))\ar[r]^{d_2}& \mathrm{H}^2(B_H,\C^*)\\
}$$

We see that a projective representation $\rho \in \mathrm{H}^1(B_H,\mathrm{PGL}_m(\C))$ comes from a linear (resp. special linear) representation if and only if $d_2\rho=0$ (resp. $d_1\rho=0$). As the diagram commutes, $d_2\rho\in i_*\mathrm{H}^2(B_H,\mathbb{Z}/m\mathbb{Z} )$, in particular it is an $m$-torsion element of $\mathrm{H}^2(B_H,\C^*)$.

Hence, if $\mathrm{H}^2(B_H,\C^*)$ has trivial $m$-torsion then $\rho$ lifts to a linear representation.
Yet, this general approach is not sharp enough for our purposes. For example  if $n>4$, $\mathrm{H}^2(F_{0,n}\Pu,\Z)$ has positive rank, see \cite{MR1217488}. As a consequence, by universal coefficient theorem, its tensor product with $\C^*$ is a subgroup of  $\mathrm{H}^2(F_{0,n}\Pu,\C^*)$ which contains torsion elements of any order.

For this reason, we use another strategy: we will lift the representation $\widehat{\mathrm{P} \rho}$ after restriction to a covering of a Zariski open set of $B_H$. We are mainly interested in finite coverings $B_H\rightarrow F_{0,n}\Pu$. Their underlying spaces are Zariski open sets of projective manifolds.

Thus, the following statement is powerful enough for our purposes. It is a known result in \'etale cohomology. We have proposed a new proof for it in \cite{CRAS}. 
\begin{theorem}\label{liftingrepr}
Let $X$ be an irreducible projective complex variety and $D$ an algebraic hypersurface in $X$. Let $\star$ be a smooth point of $X\setminus D$. 
For any representation $\rho : \pi_1(X\setminus D,\star)\rightarrow \mathrm{PSL}_m(\C)$, there exists a hypersurface $D_{\rho}$ with $D \subset D_{\rho}, \star \not \in D_{\rho}$ and a generically finite morphism $f_{\rho} : (Y_{\rho},\star_{\rho}) \rightarrow (X,\star)$ of projective varieties with basepoints, \'etale in the neighborhood of $\star_{\rho}$, such that $Y_{\rho}$ is smooth and  the pullback $$f_{\rho}^*\rho : \pi_1(Y_{\rho}\setminus f_{\rho}^{-1}(D_{\rho}),\star_{\rho})\rightarrow \mathrm{PSL}_m(\C)$$ lifts to $\mathrm{SL}_m(\C)$. 
\qed
\end{theorem}

\subsection{Conclusion}
\subsubsection{General rank}

Let us first analyze the relations between the orbit of a representation and the ones of its subrepresentations. 
\begin{lemma}\label{lemsemisimple} Fix $k$ a field.
Let $\rho=\oplus_{i\in I}~\rho_i$ be a direct sum of irreducible representations $$\rho_i: \gf \rightarrow \mathrm{GL}(W_i).$$ 
For subspaces $W_i$ of $k^m$.
Then, for any pure braid $\beta \in \pi_1(F_{0,n}\Pu,x)$, $\beta\cdot\rho=\oplus_{i\in I}~{\beta}\cdot\rho_i$. In particular, if $\beta$ stabilizes every $[\rho_i]$, it stabilizes $[\rho]$.

 Moreover, the intersection $\cap_i Stab[\rho_i]$ is a finite index subgroup of $Stab[\rho]$. In particular $[\rho]$ has finite orbit under $\pi_1(F_{0,n}\Pu,x)$ if and only if, for every $i \in I$, $[\rho_i]$ has finite orbit under $\pi_1(F_{0,n}\Pu,x)$.
\end{lemma}
\begin{proofbis}
If $\alpha_j \in \gf=\pi_1\left (\Pu \setminus \{x_1,\ldots,x_n\},x_0\right)$ is a simple loop around $x_j$, then $\alpha_j^\beta=c \alpha_j c^{-1}$, for a certain $c\in \gf$, and $(\beta\cdot\rho)(\alpha_j)=\rho(c)\rho(\alpha_j)\rho(c)^{-1}$.
For any vector $v \in W_i$, we have $\rho(c)\rho(\alpha_j)\rho(c)^{-1}\cdot v=\rho_i(c)\rho_i(\alpha_j)\rho_i(c)^{-1}\cdot v=(\beta\cdot\rho_i)(\alpha_j) v$. Using the fact that the elements $(\alpha_j)_j$ generate $\gf$ yields the first assertion of the lemma.

For any element $\beta \in Stab[\rho]$, there exists a matrix $g_{\beta}\in \mathrm{GL}_m(k)$ such that, for any $\alpha \in \gf$, $(\beta\cdot\rho)(\alpha) g_{\beta}=g_{\beta} \rho(\alpha)$. We see that, for any $i\in I$, $g_{\beta}W_i$ is an irreducible subspace of $(\beta\cdot\rho)$, in particular $[\rho_i]=[\beta\cdot \rho_k]$ for some $k\in I$. Let $J\subset I$ be such that $(\rho_{j})_{j\in J}$ is a system of representatives for the classes $([\rho_i])_{i\in I}$.
We define a morphism $Stab[\rho] \rightarrow \mathrm{Perm}(J), \beta \mapsto \sigma_\beta$ by the relation $[\beta\cdot\rho_j]=[\rho_{\sigma_{\beta}(j)}]$.
 Its kernel is a finite index subgroup of $Stab[\rho]$ that stabilizes every of the classes $[\rho_i]$.
\qed \end{proofbis}

We turn to the main result of this section.
We now consider $F_{0,n}\Pu$ as a subspace of the projective variety $(\Pu)^n$. The expression ``Zariski open" refers to open sets in the induced topology.
Recall we have defined the tautological bundle $\mathcal{T}_n$ over $F_{0,n}\Pu$ as a subset of $F_{0,n}\Pu \times \Pu$ --- see Definition $\ref{tautbundle}$. Denote $z: \mathcal{T}_n\rightarrow \Pu$ the map induced by projection to the second factor of this ambient space. Also keep in mind our notation $\gf:=\pi_1\left(\Pu\setminus \{x_1,\ldots,x_n\}, x_{0}\right)$.
\begin{theorem}\label{topthm}
Let $\rho : \gf \rightarrow \mathrm{GL}_m(\C)$ be a semisimple representation.
Suppose that  $[\rho]\in\mathrm{Hom}\left(\gf,\mathrm{GL}_m(\C)\right)/\mathrm{GL}_m(\C)$
has finite orbit under the pure braid group $\pi_1(F_{0,n}\Pu,x)$.

Then there exists a generically finite morphism $f: (U,b) \rightarrow (F_{0,n}\Pu,x)$ of smooth quasiprojective varieties, \'etale in the neighborhood of $b$,
  such that if 
\begin{itemize}
\item $p_{f} :E\rightarrow U$ is the pullback by $f$ of the bundle ${\mathcal{T}_n}\rightarrow F_{0,n}\Pu$, 
\item $\tilde{z}$ is the composition $E\rightarrow \mathcal{T}_n\stackrel{z}{\rightarrow} \Pu$,
\item $z_1 : E_b \stackrel{\simeq}{\longrightarrow} \Pu\setminus \{x_1,\ldots,x_n\}$ is the isomorphism between $E_b:=p_{f}^{-1}(b)$ and $\Pu\setminus \{x_1,\ldots,x_n\}$ induced by $\tilde{z}$, 
\item $c={z_1}^{-1}(x_{0})$,
\end{itemize}
 then there exists a representation $\hat{\rho} : \pi_1(E,c)\rightarrow \mathrm{GL}_m(\C)$ such that $\rho{z_1}_*=\hat{\rho} i_*$, where $i : E_b\rightarrow E$ is the inclusion map.
 Moreover, $\hat{\rho}$ has the same stable subspaces as $\rho$, and if certain irreducible factors of $\rho$ have trivial determinants, so can be chosen the corresponding ones for $\hat{\rho}$.
\end{theorem}
 \begin{proofbis} In the proof, we will identify $\gf$ with the fundamental groups of certain fibers of pullbacks of $\mathcal{T}_n$; these identifications are always the ones induced by the projection $z: \mathcal{T}_n\rightarrow \Pu$.
 For semidirect product structures that appear below, see section \ref{groupeloctriv}.

 Let $\rho=\oplus_{i\in I} \rho_{i}$ be the decomposition of $\rho$ in irreducible factors, with $\rho_i=\rho_{\vert W_i}$, $\C^m=\oplus_i W_i$.
 Let $H_i <\pi_1(F_{0,n}\Pu,x)$ be the stabilizer of $[\rho_i]$. By Lemma $\ref{lemsemisimple}$, $H=\cap_{i\in I} H_i$ is a finite index subgroup of $\pi_1(F_{0,n}\Pu,x)$. Thus, the covering $(B_H,b_H)\rightarrow (F_{0,n}\Pu,x)$ with fundamental group $H$  is finite.
 For this reason $B_H$ embeds as the complement of an hypersurface in a projective variety.

Consider   the bundle $E_H \rightarrow B_H$ introduced in Definition $\ref{tautbundle}$. By Corollary $\ref{extproj}$, for any $i\in I$, the projectivization of $\rho_i$ extends to a projective representation  $$\widehat{P\rho_i} :\pi_1(E_H,s_H(b))\simeq\gf\rtimes_{\varphi}\pi_1(B_H,b) \rightarrow \mathrm{PGL}(W_i).$$
Let $r_H^i :  \pi_1(B_H,b_H)  \rightarrow \mathrm{PGL}(W_i) $ be the restriction of $\widehat{\mathrm{P}\rho_i}$ to the second factor.
 By successive applications of Theorem $\ref{liftingrepr}$, there exist a Zariski open neighborhood $V$ of $b_H$ in $B_H$ and a generically finite morphism $\chi :(U,b)\rightarrow (V,b_H)$, \'etale in the neighborhood of $b$, such that if $g$ is the composition $g:U\stackrel{\chi}{\rightarrow}V\stackrel{}{\hookrightarrow} B_H$, the pullback $r_{\pi_1(U,b)}^i$ of $r_{H}^i$ by $g_*$ lifts to $\mathrm{SL}(W_i)$; for any $i\in I$.

 Let $f=f_H\circ g$ and $E \rightarrow U$ be the pullback of the bundle ${E_H}\rightarrow B_H$ by $g$; that is $E=f^*\mathcal{T}_n$. Denote $\tilde{g}: E\rightarrow {E_H}$ the natural bundle map. 
 Corollary $\ref{corlift}$ says that the existence of a lifting for $r_{\pi_1(U,b)}^i$   implies  that the pullback $\widehat{P\rho_i}_{\vert U}=\tilde{g}^*\widehat{P\rho_i}$ lifts to an extension $\widehat{\rho_i}$ of $\rho_i$
 $$\widehat{\rho_i} : \pi_1\left(E,s_{f}(b)\right)\simeq \gf\rtimes_{\varphi}\pi_1(U,b)\rightarrow \mathrm{GL}(W_i);$$
 
  with values in $\mathrm{SL}(W_i)$ if $\mathrm{Im}(\rho_i)\subset \mathrm{SL}(W_i)$. We can then define an extension of $\rho$ as $\hat{\rho}:=\oplus_i\widehat{\rho_i}$.
  The map $f$ is the one announced in the statement of the theorem. \qed \end{proofbis}

\subsubsection{Rank two representations}
In view of our application to Garnier systems, we give a specific enhanced version of Theorem $\ref{topthm}$ for rank two representations.
\begin{theorem}\label{topthm rang2}
Suppose that $[\rho]\in\mathrm{Hom}\left(\gf,\mathrm{GL}_2(\C)\right)/\mathrm{GL}_2(\C)$
has finite orbit under the pure braid group $\pi_1(F_{0,n}\Pu,x)$.

Then there exists a generically finite morphism $f: (U,b) \rightarrow (F_{0,n}\Pu,x)$ of smooth quasiprojective varieties, \'etale in the neighborhood of $b$,
  such that if 
\begin{itemize}
\item $p_{f} :E\rightarrow U$ is the pullback by $f$ of the bundle ${\mathcal{T}_n}\rightarrow F_{0,n}\Pu$, 
\item $\tilde{z}$ is the composition $E\rightarrow \mathcal{T}_n\stackrel{z}{\rightarrow} \Pu$,
\item $z_1 : E_b \stackrel{\simeq}{\longrightarrow} \Pu\setminus \{x_1,\ldots,x_n\}$ is the isomorphism between $E_b:=p_{f}^{-1}(b)$ and $\Pu\setminus \{x_1,\ldots,x_n\}$ induced by $\tilde{z}$, 
\item $c={z_1}^{-1}(x_{0})$,
\end{itemize}
 then there exists a representation $\hat{\rho} : \pi_1(E,c)\rightarrow \mathrm{GL}_2(\C)$ such that $\rho{z_1}_*=\hat{\rho} i_*$, where $i : E_b\rightarrow E$ is the inclusion map.
  
  If $\rho$ takes values in $\mathrm{SL}_2(\C)$ so can be chosen $\hat{\rho}$.
\end{theorem}
\begin{proofbis}
For irreducible and abelian  non scalar representations, this is direct application of Theorem $\ref{topthm}$ and Lemma $\ref{lemsemisimple}$.

 Any scalar representation $\rho: \gf \rightarrow \C^*$ is fixed by the action of $\pi_1(F_{0,n}\Pu,x)$, because the loops around the punctures $(x_i)$ generate the fundamental group $\gf$. We can thus extend any such $\rho$ to $\hat{\rho}: \gf\rtimes \pi_1(F_{0,n}\Pu,x)\rightarrow \C^*$ by the formula $\alpha \beta \mapsto \rho(\alpha)$, the map $f$ can be taken to be the identity map of $F_{0,n}\Pu$.

We now treat the reducible non abelian case.
After conjugation and tensor product with a scalar representation, we may suppose $\rho$ takes its values in the subgroup
$\mathrm{Aff}(\C):=\{\left [\begin{smallmatrix}\lambda&\tau\\0&1\end{smallmatrix}\right] , \lambda \in \C^*, \tau \in \C \}$ of $\mathrm{GL}_2(\C)$. In this case, we easily see that for any $\beta \in H=Stab[\rho]$ there exists a matrix $A_{\beta}\in \mathrm{Aff}(\C)$ satisfying $(\beta\cdot\rho)A_{\beta}=A_{\beta} \rho$. 

Of course we will interpret $\mathrm{Im}(\rho)<\mathrm{Aff}(\C)$ as a group of transformations of the affine complex line $\C$. By non abelianity,  this group does not fix any point in $\C$ and is not a group of translations. This means that we have two non trivial elements $\rho(\alpha_i), i=1,2$ of $\mathrm{Im}(\rho)$
fixing two distinct points in $\C$. The  element $A_{\beta}$ is thus characterized as the unique element of $\mathrm{Aff}(\C)$ sending the fixed point of $\rho(\alpha_i)$ to the fixed point of $(\beta\cdot\rho)(\alpha_i)$ for $i=1,2$.
 From uniqueness of $A_\beta$, we see that $\beta\mapsto A_{\beta}, H \rightarrow \mathrm{SL}_2(\C)$ is a homomorphism. We then extend $\rho$ to $\gf \rtimes_{\zeta} H$ by the formula $\alpha \beta \mapsto \rho(\alpha) A_{\beta}$. The map $f$ is the map $f_H: (B_H,b_H)\rightarrow (F_{0,n}\Pu,x)$.
\qed \end{proofbis}
\subsubsection{Non vacuity}\label{non vacuity}
Just to show our construction is not empty, we give two basic examples of finite orbits.
First, notice that any representation $\Lambda_n \rightarrow \mathrm{GL}_m(\C)$ with finite image yields a finite orbit. For any $m>0$, the two generator group $\mathfrak{S}_{m+1}$ has an irreducible representation of rank $m$. Hence for any $n\geq3$, $m>0$, we have a finite orbit arising from an irreducible representation.

One of the referees made the following interesting remark. As the total space $\mathcal T_n$ of the tautological bundle over $F_{0,n}\Pu$ is nothing but $F_{0,n+1}\Pu$, for any $\rho: \Lambda_n\rtimes \pi_1(F_{0,n}\Pu) \simeq \pi_1(F_{0,n+1}\Pu) \rightarrow \mathrm{GL}_m(\C)$, the  induced representation of $\Lambda_n$ yields a size $1$ orbit under $\pi_1(F_{0,n}\Pu)$.
\section{Riemann-Hilbert correspondence}\label{secRH}
 For an introduction to flat connections, we refer to \cite[sec. $17$]{MR2363178}.\\
Let $X$ be a complex manifold. We denote $\mathcal{O}_X$, $\mathcal{M}_X$ the sheaves of holomorphic and  meromorphic functions on $X$, respectively.
A rank $m$ meromorphic connection \emph{over} $X$ is map of $\C$-vector space sheaves $\nabla : {\bf{V}}\rightarrow {\mathcal{M}_X\otimes \bf{V}}$, such that $\bf{V}$ is the sheaf of sections of a rank $m$ holomorphic vector bundle $V$ over $X$ and  for any $x\in X$, $f\in{\mathcal O}_{X,x}$, $s\in {\bf{V}}_x$,  $\nabla(f\cdot s)=f\cdot \nabla(s)+ds \otimes f$.

Thus, the data of $V$ is contained in the one of $\nabla$. However, in need of explicit reference to the vector bundle, we shall use the notation $(V,\nabla)$ and say that $\nabla$ is a connection \emph{on} $V$ or that $V$ is endowed with the connection  $\nabla$.

In \cite{MR0417174}, for a complex manifold $X$ containing a normal crossing hypersurface $D$, Deligne realizes any represention $\rho: \pi_1(X\setminus D)\rightarrow \mathrm{GL}_m(\C)$ as the monodromy representation of a logarithmic flat connection $\nabla_{\rho}$ over $X$ with polar locus $D$.
 He proceeds \textit{via} a patching procedure, glueing some local models that satisfy a non resonance property. We want to describe a slight modification of his construction that allows a wider range of local models and gives some freedom to prescribe the various transversal types for $\nabla_{\rho}$, see Definition $\ref{deftransversaltype}$ below.
This will be needed to prove  Theorem $\ref{maintheorem}$.

 A complete description of all possible choices of logarithmic extensions would be interesting. It would certainly involve the topology of the polar divisor.

\subsection{Local theory of logarithmic flat connections}
\subsubsection{First results} \label{firstlocalresults}
We must first understand the local behavior of a logarithmic connection at smooth points of the polar divisor. The following allows reduction of this study to a one dimensional base space. 
\begin{proposition} \label{local constancy}
Let $w_1,\ldots, w_d$ be the coordinates of $\C^d$. Let $U$ be neighborhood  of $0\in \C^d$.
Let $\nabla$ be a flat logarithmic connection on $\mathcal{O}^{\oplus m}_U$,  with polar divisor $w_1=0$. Then up to reduction of $U$, $U$ is a product $\D_r\times V$ for a disk $\D_r=\{w_1\in \C, \vert w_1\vert<r\}$ and 
after a holomorphic gauge transformation, $\nabla$ is  the pullback of a connection $\nabla_0$ on $\mathcal{O}^{\oplus m}_{\D_r}$ by the projection $w\mapsto w_1$.
\begin{proofbis}
This is a special case of \cite[Theorem $5$]{MR0436189}.
\qed \end{proofbis}
\end{proposition}
On $(\C,0)$ we use the standard coordinate of $\C$, which we denote $w$.
 We denote $\Re(\lambda)$ the real part of a complex number $\lambda$.
\begin{definition}\label{reduced}
We say that a logarithmic connection $Y\mapsto dY-A(w)\frac{dw}{w}Y$ on $\mathcal{O}^{\oplus m}_{(\C,0)}$, is reduced if 
$A$ has a block diagonal shape $$A(w)=\begin{bmatrix}
                    B^1(w)& 0 &\cdots&0 \\
                    0 & \ddots &\ddots& \vdots \\
                 \vdots& \ddots &  \ddots&0\\
                     0& \cdots &0&B^r(w)
                  \end{bmatrix} $$
                  such that the blocks $B^i(w)=\sum_{k\geq 0} B^{i}_k w^k$ satisfy the following conditions.
\begin{enumerate}
\item The blocks $B^i(w)$ are upper triangular and the matrices $\Gamma^i:=B^{i}_0$ are in Jordan form,
\item if $i\neq j$, for any eigenvalue $\lambda$ of $\Gamma^i$ and any eigenvalue $\mu$ of $\Gamma^j$, $\lambda-\mu\not \in \mathbb{Z}$.
\item For any pair of eingenvalues $\lambda,\mu$ of $\Gamma^i$, $\lambda-\mu\in \mathbb{Z}$.
\item The eigenvalues of $\Gamma^{i}$ are in decreasing order on the diagonal: $\Re({\Gamma^i_{l,l}})$ is a (not necessarily strictly) decreasing sequence.
\item \label{monome} The $(u,v)$ entry of $B^i$ is monomial, of the form $cw^k$; with $k=\Gamma^i_{u,u}-\Gamma^i_{v,v}$ and $c\in \C$.
\end{enumerate}
\end{definition}
Notice condition $(\ref{monome})$ of the definition implies $w^{\Lambda^i} B_k^i w^{-\Lambda^i}=B_k^iw^k$, where $\Lambda^i$ is the diagonal part of $\Gamma^i$. If $\Lambda$ is the diagonal part of $A(w)=\sum_k A_k w^k$, this yields
$w^{\Lambda} A_k w^{-\Lambda}=A_k w^k.$ We also remark that $A(w)$ is a polynomial and the expression $A(1)$ makes sense.
\begin{theorem}[Poincar\'e-Dulac-Levelt]\label{thmPDL}
Any germ of logarithmic connection on $\mathcal{O}^{\oplus m}_{(\C,0)}$ is holomorphically gauge equivalent to  a reduced connection.
\end{theorem}
\begin{proofbis}
See \cite[Sec. 1$6C-16D$]{MR2363178}.
\qed \end{proofbis}

For a real number $a$, $\lfloor a \rfloor$ denotes the biggest integer $q$ with $a\geq q$.
\begin{proposition}\label{birat}
Let $\nabla :Y\mapsto dY- A(w)\frac{dw}{w}Y$ be a reduced logarithmic connection on $\mathcal{O}^{\oplus m}_{(\C,0)}$. Let  $\Lambda$ be the  diagonal part of $A$. Let $L$ be the diagonal matrix defined by ${L}_{u,u}=\lfloor\Re({\Lambda}_{u,u}) \rfloor$. The meromorphic gauge transform $Y=w^{L} Z$ transforms this connection into a connection $\nabla^0 : Z \mapsto dZ-C\frac{dw}{w}Z$ with constant $C$ given by $C=A(1)-L$. In particular, any eigenvalue $\lambda$ of $C$ satisfies $\Re(\lambda)\in [0,1)$.
\begin{proofbis}
This computation is already mentioned in \cite[\S 4]{MR0436189}. If $Y$ is a horizontal section of $\nabla$ and $Z=w^{-L}Y$  we have $$dZ=d(w^{-L}Y)= (-L+w^{-L}A(w)w^L)\frac{dw}{w} Z.$$
If $\Lambda$ is the diagonal part of $A(w)$, the block diagonal structure of $A(w)$ allows to see  $w^{-L}A(w)w^L=w^{-\Lambda}A(w)w^{\Lambda}$. We mentioned above that $A(w)=\sum_k A_k w^k$ can be rewritten $A(w)=\sum_k w^{\Lambda}A_k w^{-\Lambda}$. Hence $w^{-\Lambda}A(w)w^{\Lambda}=\sum_k A_k=A(1)$.
\qed \end{proofbis}
\end{proposition}
\begin{definition} 
We define two spaces of connections on $\mathcal{O}^{\oplus m}_{(\C,0)}$. \begin{itemize}
\item The space of nonresonant Euler connections:
$$Eul^0:=\cup_m\{Z \mapsto dZ- C\frac{dw}{w}Z, C\in \mathrm{Mat}_m(\C), \lambda,\mu \in \mathfrak{S}(C) \Rightarrow \lambda-\mu \not \in \mathbb{Z}^* \},$$
where $\mathfrak{S}(C)$ is the spectrum of the matrix $C$.
\item The space $Red$ of reduced logarithmic connections. 
\end{itemize}
Proposition $\ref{birat}$ defines a map
$$\begin{array}{lccc} eul : & Red &\longrightarrow& Eul^0\\
							&\nabla&\longmapsto&\nabla^0
\end{array} $$
The effect of $eul$ on connection matrices shall be denoted $A\mapsto \varepsilon(A)$.
\end{definition}
Notice $\varepsilon(A)$ has the same block diagonal structure as $A$, with triangular blocks $C^i=B^i(1)-L^i$; $L^i$ diagonal with ${L^i}_{u,u}=\lfloor\Re({\Lambda^i}_{u,u}) \rfloor$.
\subsubsection{Local automorphisms}\label{localaut}
We want to investigate the automorphisms of a reduced connection in relation with the ones of its corresponding Euler connection.

\begin{proposition}\label{Euler case}
Let $\nabla^0 :Z \mapsto dZ- C\frac{dw}{w}Z$ be an element of $Eul^0$.
Any bimeromorphic automorphism $\phi$ of $\mathcal{O}^{\oplus m}_{(\C,0)}$ that preserves $\nabla^0$ ($\phi^*\nabla^0=\nabla^0$) has the form $Z\mapsto g \cdot Z$, for a constant matrix $g$.

\begin{proofbis}
This is inspired by \cite[Lemme $3$]{MR2077648}.
Such an automorphism has the form $Z\mapsto G(w) \cdot Z$, with $G=\sum_{k=r} G_k w^k$ a meromorphic function with values in $\mathrm{GL}_m(\C)$.
the relation $\phi^*\nabla^0=\nabla^0$ reads $$[G,C]dw+wdG=0$$ or
$$\sum_{k=r} ([G_k,C]+k G_k)w^k=0$$
Each term of the left hand side must be zero and we have, for $k\geq r$, $G_k C =(C-k \mathrm{Id})G_k$. Such identities imply that, for $k\geq r$, if $G_k\neq 0$, $C$ and $(C-k \mathrm{Id})$ have a common eigenvalue $\lambda=\mu-k$.
Nonresonance of $C$ shows $G_k=0$ for $k\neq 0$. Hence, $G=G_0$ is constant. 
\qed \end{proofbis}
\end{proposition}
We immediately derive the following, which  is contained in \cite[Lemma $1$]{MR2113080}.
\begin{corollary} \label{lemviktor}
Let $\nabla :Y\mapsto dY- A(w)\frac{dw}{w}Y$ be a reduced logarithmic connection on $\mathcal{O}^{\oplus m}_{(\C,0)}$.
Let $\Lambda$ be the diagonal part of $A$, then a bimeromorphic automorphism $Y\mapsto G(w)\cdot Y$ of $\mathcal{O}^{\oplus m}_{(\C,0)}$ preserves $\nabla$ if and only if  $G(w)=w^{\Lambda} g w^{-\Lambda}$ for $g$ a constant matrix satisfying $[g,\varepsilon(A)]=0$. \qed
\end{corollary}
Suppose $A$ is as in Definition $\ref{reduced}$.
 Any matrix $g$ satisfying $[g,\varepsilon(A)]=0$ must fix the characteristic spaces of $\varepsilon(A)$. We thus have the same block diagonal structure for $g$, $\varepsilon(A)$ and $A$. Hence, a matrix $g$ satisfies $[g,\varepsilon(A)]=0$, if and only if $g$ has  the same block diagonal shape as $A$
  $$g=\begin{bmatrix}
                    g^1& 0 &\cdots&0 \\
                    0 & \ddots &\ddots& \vdots \\
                 \vdots& \ddots &  \ddots&0\\
                     0& \cdots &0&g^r
                  \end{bmatrix} $$
                  and, for every $i$, $[g^i,\varepsilon(B^i)]=0$.

The holomorphicity of the matrix $G(w)=w^{\Lambda} g w^{-\Lambda}$ is then
tantamount to 
 $$(\star):~~~g_{u,v}=0\mbox{ for every }(u,v)\mbox{ such that }\Lambda_{v,v}-\Lambda_{u,u} \in \mathbb{N}^*.$$

For instance, if $g$ is upper triangular then $G(w)$ is holomorphic with holomorphic inverse.

\subsection{Transversal types}

\begin{prop/def}\label{deftransversaltype}
Let $X$ be a complex manifold and $D$ be an analytic hypersurface in $X$. Suppose we have a flat logarithmic connection  $\nabla$ on a vector bundle over $X$, with polar locus $D$. 
Let $D_0$ be an irreducible component of $D$. Let $\gamma : \D\rightarrow X$ be an embedding of the unit disk in $X$, transverse to $D$, such that $\gamma(0)$ is a smooth point of $D_0$. 
Then the gauge isomorphism class of the germ of $\gamma^*\nabla$ at $0$ is independent  of $\gamma$. 

\textnormal{This isomorphism class is called the \emph{transversal type} of $\nabla$ on $D_0$.}
\end{prop/def}
\begin{proofbis}
Let $p$ be a smooth point of $D_0$, Proposition $\ref{local constancy}$ allows to find coordinates $(w_1,\ldots,w_d): U\rightarrow \D\times V,$ $w_1\in \D$, $(w_2,\ldots,w_d) \in V$ centered at $p$ and such that $\nabla_{\vert U}$ is, up to gauge isomorphism, the pullback by $w_1$ of a logarithmic connection $\tilde{\nabla}$ on $\D$ in reduced form; $\tilde{\nabla} : Z \mapsto dZ-A(w_1) \frac{dw_1}{w_1} Z$.
By connectedness of $D_0$, it is enough to prove the result for curves $\gamma$ with $\gamma(0)\in U$.

 Let $\gamma : \D \rightarrow X$ be as prescribed above, with $\gamma(\D)\in U$. By transversality, we have $w_1\circ\gamma(t)=t u(t)$ with $u(0)\neq 0$ and, up to gauge isomorphism, the germ of $\gamma^*\nabla$ at $t=0$ is
$Z \mapsto dZ-A(tu) \dlog(tu) Z$. This connection is gauge isomorphic to $ \tilde Z\mapsto d \tilde Z-A(t)\frac{dt}{t} \tilde Z$. Indeed,  let $L$ and $C$ be the matrices defined from $A$ as in Proposition $\ref{birat}$; it is readily seen that the meromorphic gauge transform $\tilde Z=t^Lu^{-C}(ut)^{-L} Z$ relates both connections. Yet, we have to check it is a holomorphic gauge isomorphism.
To see this, rewrite $t^Lu^{-C}(ut)^{-L}=t^Lu^{-C}t^{-L}u^{-L}$, and note $t^Lu^{-C}t^{-L}=t^{\Lambda}u^{-C}t^{-\Lambda}$ is holomorphic with holomorphic inverse, because $u^{-C}$ is upper triangular, with the same block diagonal structure as $A$.
\qed \end{proofbis}

\begin{definition}
Let $\nabla :Y\mapsto dY- A(w)\frac{dw}{w}Y$ be a logarithmic connection on $\mathcal{O}^{\oplus m}_{(\C,0)}$. We will say that $\nabla$ is a \emph{mild transversal model} if any meromorphic gauge transform $Y\mapsto G(w)\cdot Y$ of $\mathcal{O}^{\oplus m}_{(\C,0)}$ which preserves $\nabla$ is holomorphic.
The gauge isomorphism class of $\nabla$ is then called a \emph{mild transversal type}.
\end{definition}
Notice that if $\nabla$ is a mild transversal model and $G(w)$ preserves $\nabla$, then so do $G^{-1}(w)$ and $G^{-1}(w)$ must also be holomorphic.
Also, a logarithmic connection on $\mathcal{O}^{\oplus m}_{(\C,0)}$ is a mild transversal model if and only if it defines a mild transversal type.
Thus, we only need to investigate mildness for reduced logarithmic connections on $\mathcal{O}^{\oplus m}_{(\C,0)}$. The mildness for reduced connections has been described in the end of section $\ref{firstlocalresults}$ (Corollary $\ref{lemviktor}$ and subsequent comments).
\begin{example} \label{nonmild}
Fix $n\in \mathbb{Z}_{>0}$, $\tau\in \C$. The reduced logarithmic connection $Z\mapsto dZ-\left [ \begin{smallmatrix}n/2+\tau&0\\0&-n/2+\tau\\ \end{smallmatrix}\right]  \frac{dw}{w}Z$ on $\mathcal{O}^{\oplus 2}_{(\C,0)}$ does not define a mild transversal type. Indeed, the gauge transform given by $G(w)=\left[\begin{smallmatrix}0&w^n\\w^{-n}&0\\ \end{smallmatrix}\right]$ preserves $\nabla$.
\end{example} We now present two families of such mild transversal models. Example $\ref{nonmild}$ shows that their rank two specializations cover all rank two mild tranversal types.

\begin{lemma}
Any element of $Eul^0$ is a mild transversal model.
\begin{proofbis}
This is a trivial consequence of Proposition $\ref{Euler case}$.
\qed \end{proofbis}
\end{lemma}
Maybe a more original example is the following.

\begin{lemma}\label{original}
Let $\nabla : Y\mapsto dY -A(w)\frac{dw}{w}Y$ be a logarithmic reduced connection on $\mathcal{O}^{\oplus m}_{(\C,0)}$.
If, for every $0<u<m$, $A_{u,u+1}(1)\neq 0$ then $\nabla$ is a mild  transversal model.
\begin{proofbis}
Let $\Lambda$ be the diagonal part of $A$ and $C:=\varepsilon(A)$.
Let $g$ be a constant matrix of size $m$ satisfying $[g,C]=0$.
The matrix $C$ is upper triangular, thus
the relation $[g,C]_{u,v}=0$ reads $$r_{u,v}:~~~~~~\sum_{l=1}^{v-1} g_{u,l}C_{l,v}=\sum_{l=u+1}^{m} C_{u,l}g_{l,v}$$
We can proceed by induction to show $g$ is upper triangular; yielding the required holomorphy.
The induction hypothesis is that, up to column $k$ included, $g$ is like a triangular matrix, namely
$$\mathcal{H}(k):\mbox{ for every }j\mbox{ such that } 0<j\leq k,\mbox{ for every }i>j,   g_{i,j}=0.$$
We have trivially $\mathcal{H}(0)$; let us prove $\mathcal{H}(k)$ implies $\mathcal{H}(k+1)$, provided $k<m$.

By the induction hypothesis the relation $r_{u,k+1}$ reduces to $$\sum_{l=u}^{k} g_{u,l}C_{l,v}=\sum_{l=u+1}^{m} C_{u,l}g_{l,k+1},$$
In particular, for $u> k$, the right hand side of this equation is zero. The family of equations $(r_{u,k+1})_{m>u> k}$ is a triangular homogeneous linear system in the unknowns $(g_{l,k+1})_{m\geq l>k+1}$.
The diagonal coefficient of line $u$ is $C_{u,u+1}=A_{u,u+1}(1)$ and is nonzero by assumption, this proves $g_{l,k+1}=0$ for $m\geq l>k+1$, which is exactly $\mathcal{H}(k+1)$.
The proof is complete. 
\qed \end{proofbis}
\end{lemma}
Notice, that the assumptions of Lemma $\ref{original}$ actually imply that $A$ has only one block $B^1$ in the decomposition of Definition $\ref{reduced}$.

If $A$ has several blocks  $B^i$ such that the connections $\nabla_i: Z_i \mapsto dZ_i-B^i(w)\frac{dw}{w}Z_i$ are all mild transversal models, then $\nabla=\oplus_i \nabla_i$ is also a mild transversal model. Indeed, we have seen in section $\ref{localaut}$ that any bimeromorphic automorphism of $\nabla$ has the same block diagonal structure as $A$.

\subsection{Logarithmic Riemann-Hilbert}
We start with a local result that requires no mildness assumption.
\begin{lemma}\label{localRH}
Let $w_1,\ldots,w_d$ be the standard coordinates of $\C^d$.\\
Let $\Delta=\{(w_1,\ldots,w_d)\in \C^d\vert~|w_j|<2, j=1,\ldots,d\}$. Let $D_j=\{w_j=0\}$ and $D=\cup D_j$.
Denote $b=(1,\ldots,1)$.
Let $\rho : \pi_1(\Delta\setminus D,b)\rightarrow  \mathrm{GL}_m(\C)$ be an antirepresentation.
Let $t \mapsto \alpha_j(t)$ be the  loop with coordinate functions $w_j(t)=exp(2i \pi t); w_k(t)\equiv 1, j\neq k.$ 
Choose logarithmic reduced connections $\nabla_j$ with monodromy along $w_j(t)=exp(2i \pi t)$ conjugate to $\rho(\alpha_j)$.

Then there exists a logarithmic connection $\nabla$ on $\mathcal{O}^{\oplus m}_U$ with polar locus $D$, monodromy equal to $\rho$ and transversal type on $D_j$ given by $\nabla_j$.
\begin{proofbis}
We apply Proposition $\ref{birat}$ to each of the reduced logarithmic connections $\nabla_j:Z\mapsto dZ-A_j(w_i) \frac{dw_j}{w_j}Z$.
Let $C_j:=\varepsilon(A_j)$ and $L_j:=A_j(1)-C_j$.

 The matrix $w_j^{L_j}$ is the identity in restriction to $w_j=1$, thus we have exactly (with no conjugation) the same monodromy along $w_j(t)=exp(2i \pi t)$ for $\nabla_j$ and $eul(\nabla_j)$, namely $exp(2i\pi C_j)$.
By hypothesis, there exists $G_j$ such that $G_j exp(2i \pi C_j) G_j^{-1}=\rho(\alpha_j)$. Moreover, $C_j$ has real parts of eigenvalues contained in $[0,1)$. Therefore, by \cite[lemme $2$]{MR2077648}, the matrix $R_j=G_jC_jG_j^{-1}$ is the unique matrix with real parts of eigenvalues contained in $[0,1)$ and $exp(2i\pi R_j)=\rho(\alpha_j)$. Then, by the same lemma  and commutativity of the $\alpha_j$'s, the matrices $R_j$ commute.

 Hence, the connection $ Z \mapsto dZ-\sum_j R_j\frac{dw_j}{w_j}Z$ is flat and its monodromy  is exactly $\rho$.
We transform this connection by $Y=\prod_j (G_j w_j^{L_j}G_j^{-1}) Z$ to obtain a new flat logarithmic connection $\nabla : Y\mapsto dY-\sum B_j(w)\frac{dw_j}{w_j}Y$. We see that,  in restriction to $w_k=1, k\geq 2 $, this transformation is $Y=(G_1 w_1^{L_1}G_1^{-1}) Z$ and the restriction of $\nabla$
is $Y \mapsto dY-G_1 A_1(w_1)\frac{dw_1}{w_1} G_1^{-1}Y$, which is holomorphically equivalent to $\nabla_1$. Obviously, we can make the same reasoning for other polar components; $\nabla$ is the sought connection.
\qed \end{proofbis}
\end{lemma}

\begin{definition}
Let $X$ be a complex manifold, let $D \subset X$ be a normal crossing hypersurface. Let $\rho: \pi_1(X\setminus D,b) \rightarrow \mathrm{GL}_m(\C)$ be a representation.
Let $\nabla$ be a flat rank $m$ connection over $X\setminus D$ with monodromy $\rho$.
Let $D_i$ be an irreducible component of $D$ and let $\nabla_i$ be a rank $m$ connection over the germ of disk $(\D,0)$. Let $\gamma : (\D,0) \rightarrow X$ be a germ of curve cutting $D_i$ transversely at a smooth point of $D$.
We say that $\nabla_i$ is \emph{compatible with $\rho$} if $\nabla_i$ and $\gamma^*\nabla$ have the same monodromy representation up to conjugation.
\end{definition}

\begin{theorem}\label{RHthm}
Let $X$ be a complex manifold, let $D \subset X$ be a normal crossing hypersurface and let $$\rho: \pi_1(X\setminus D,b)\rightarrow  \mathrm{GL}_m(\C)$$ be an antirepresentation. 
Denote $(D_i)_{i \in I}$ the irreducible components of $D$ and choose, for every $i\in I$, a mild transversal model $\nabla_i$ compatible with $\rho$.

Then there exists a flat logarithmic connection over $X$ with polar locus $D$, transversal type on $D_i$ given by $\nabla_i$ and monodromy representation  given by $\rho$.
This connection is unique up to vector bundle isomorphism.

\begin{proofbis} We first prove the existence part.
We shall give connections on open subspaces of $X$ which cover $X$. Then we will show that they glue together in a logarithmic flat connection over $X$.

By suspension, we have a holomorphic flat connection $\nabla_0$ on $U_0=X\setminus D$ with monodromy $\rho$.
In the neighborhood $U_p$ of any smooth point $p$ of $D$ we have a flat connection $\nabla_p$ given by the mild transversal model $\nabla_i$, $p\in D_i$. In the neighborhood $U_p$ of a non smooth point $p\in D$, Lemma $\ref{localRH}$ gives us a flat connection $\nabla_p$ with transversal type on $D_i$ given by $(\nabla_i)$ and with the same monodromy as ${\nabla_0}$ on $ U_p\cap U_0$.

For any $p \in D$ we have, by coincidence of monodromy, an isomorphism $$\phi_{0,p}^*\nabla_{0\vert U_0 \cap U_p}=\nabla_{p\vert U_0\cap U_p}.$$

For any pair of points $q,p\in D$ such that $U_p \cap U_q\neq \varnothing$, we want to define isomorphisms $\phi_{q,p}^*\nabla_{q\vert U_q \cap U_p}=\nabla_{p\vert U_0\cap U_p}$. Outside $D$ such an isomorphism is given by $\phi_{q,p}=\phi_{0,q}^{-1}\circ \phi_{0,p}$. Let us show it extends holomorphically to $D\cap U_p \cap U_q$: by regularity $\phi_{q,p}$ is meromorphic at $D$, \cite{MR0417174}. Thus, it is enough to check its restriction to any curve crossing $D$ transversely at a general smooth point of $D\cap U_p \cap U_q$ is holomorphic, because the indeterminacy set of a meromorphic function is a proper analytic subset of its polar locus. Such a verification is automatic by the mildness of the transversal models $\nabla_i$.

We thus have a set of holomorphically invertible gluing maps $\phi_{\alpha,\beta}$ defined on $U_{\alpha} \cap U_{\beta}$, which are compatible in the sense: $$\phi_{\alpha,\gamma}=\phi_{\alpha,\beta} \circ \phi_{\beta,\gamma}~;~~~ \phi_{\alpha,\beta}=\phi_{\beta,\alpha}^{-1}.$$
Glueing the vector bundles with connections $(\nabla_{\alpha})$ using these maps, we obtain the desired flat logarithmic connection over $X$.

The uniqueness assertion is obtained as follows:  for two connections satisfying the conclusion of the theorem, by coincidence of monodromy, we have an isomorphism outside $D$; it extends holomorphically to $D$ by mildness of the $\nabla_i$'s.
\qed \end{proofbis}
 \end{theorem}\section{Isomonodromic deformations}\label{isomdef}

We introduce some terminology on isomonodromic deformations and then prove our main result.
\subsection{Definition}

Let $T$ be a complex manifold and consider a flat rank $m$ connection $\nabla$ over $T\times \Pu$ with polar locus given by holomorphic sections of $T\times \Pu\rightarrow M$, with disjoint images.
For $t\in T$, consider the restriction $\nabla_t=\nabla_{\vert \{t\}\times \Pu}$. 
 For $t^0\in T$, we can obviously see the family $(\nabla_t)_{t\in T}$ as a deformation of the connection $\nabla_{t^0}$ over $\Pu$.
 By flatness of $\nabla$, we see that, for a given loop $\alpha$ in the complement of the poles of $\nabla_{t^0}$ in $\Pu$, the monodromy of $\nabla_{t_0}$ over $\alpha$ is the same as the one of 
 $\nabla_{t^1}$ if $t^1$ is close enough to $t^0$ in $T$.
 
 For this reason, we call the family $(\nabla_t)_{t\in T}$ an \emph{isomonodromic deformation}.
 Conversely, we will use this wording \emph{exclusively} in the above situation, where the family $\nabla_t$ is obtained by restriction of a logarithmic connection $\nabla$.


\subsection{Universal deformation}
Denote \[\begin{array}{rccl}r:& (B_0,b_0)&\rightarrow& (F_{0,n}\Pu,x)\\
						&\tilde{y}&\mapsto&(y_1,\ldots,y_n) \end{array}\] the universal cover of $F_{0,n}\Pu$.
In Definition $\ref{tautbundle}$, we considered the tautological fiber bundle $E_0$ obtained by restriction of
$$\begin{array}{lccc} p:& B_0\times\Pu& \rightarrow &B_0\\
							& (\tilde{y},z)&\mapsto&\tilde{y}
\end{array}$$ to the open set $(r\times id_{\Pu})^{-1}(\cap_i\{y_i\neq z\})$.

For $\nabla$ a logarithmic connection on a rank $m$ holomorphic vector bundle $V$ over $\Pu$ with polar locus $\{x_1,\ldots,x_n\}$, Malgrange  \cite{MR0728431} has constructed a  rank $m$ vector bundle $\hat{V}$ over $B_0\times \Pu$ endowed with a logarithmic flat connection $\hat{\nabla}$ with poles at $\{y_i=z\}$
and such that $(V,\nabla)=i^*(\hat{V},\hat{\nabla})$, for $i : \Pu \rightarrow B_0\times \Pu, z\mapsto (b_0,z)$. The pair $(\hat{V},\hat{\nabla})$ is unique up to bundle isomorphisms. Actually, if $\nabla$ has only mild transversal types, this existence and uniqueness statement is a direct consequence of Theorem $\ref{RHthm}$.

 This connection $\hat{\nabla}$ has the following universal property.
 \begin{proposition}\label{propuniv}
Let  $B$ be a simply connected complex manifold  and $f: (B,b) \rightarrow (F_{0,n},x)$.
Let $$D= \cup_i \{\left((y_1,\ldots,y_n),z\right)\in \left (\Pu\right)^n\times \Pu\vert y_i=z\}.$$
Let $\tilde{f}$ be the lift of $f$ to $B_0$.
Let $\nabla$ be a logarithmic connection over $\Pu$, with polar locus $\{x_1,\ldots,x_n\}$ and define $j : \Pu \rightarrow B \times \Pu, z\mapsto(b,z)$.

Let $\tilde \nabla$ be a flat logarithmic connection over $B \times \Pu$ with polar locus $(f\times~id_{\Pu})^{-1}D$ and such that $j^*(\tilde{\nabla})=\nabla$. 

Then $\tilde \nabla$ is isomorphic to the pullback $(\tilde{f}\times id_{\Pu})^*\hat{\nabla}$.
\begin{proofbis}
See \cite{MR0728431}; again if $\nabla$ has only mild transversal types, this is an easy corollary of Theorem $\ref{RHthm}$.
\qed \end{proofbis}
 
 \end{proposition}
 Because of this property, $\hat{\nabla}$ is called \emph{the universal isomonodromic deformation of~$\nabla$}.
 \subsection{Main result}
 \begin{definition}\label{defalg}
Let $\nabla$ be a logarithmic connection on a rank $m$ holomorphic vector bundle over $\Pu$.
We will say that \emph{the germ of universal isomonodromic deformation of $\nabla$ is algebraizable} if there exist: 
\begin{enumerate}
\item a birational morphism $\pi: X\rightarrow Y\times \Pu$ between projective complex manifolds,
\item a neighborhood $\Delta$ of $b_0$ in $B_0$ and an open holomorphic embedding $$\psi : \Delta  \rightarrow Y, $$ whose image does not intersect the projection in $Y$ of the indeterminacy locus of $\pi^{-1}$. 
\item a holomorphic embedding $\Psi: \Delta\times \Pu \rightarrow X$ such that the following diagram commutes, $q : X \rightarrow Y$ being the ruling of $X$.
 $$\xymatrix{
 \Delta \times \Pu \ar[r]_{\Psi} \ar[d]_{p}&X \ar[d]^{q}\\
    \Delta \ar[r]^{\psi} & Y}$$

\item \label{algext} a rank $m$ vector bundle  $\tilde{V}$ on $X$ endowed with a meromorphic flat connection $\tilde{\nabla}$,
\item a connection isomorphism $(\hat{V},\hat{\nabla})_{\vert \Delta\times \Pu} \simeq \Psi^*(\tilde{V},\tilde{\nabla})$,
\end{enumerate}
If the connection $\tilde{\nabla}$ of $\ref{algext}$. can be chosen logarithmic, we will say that  that the germ of universal isomonodromic deformation of $\nabla$ is \emph{logarithmically algebraizable}.
\end{definition}
\begin{remark}
In this definition we refer to the idea of algebraicity for the following reasons.
\begin{enumerate}
\item Any holomorphic vector bundle $\tilde{V}$ on a projective manifold  $X$ comes from an algebraic one. This means we can find a covering $(U_i)$ of $X$ by Zariski open sets and holomorphic trivializations of $\tilde{V}$ relative to this covering in order to have transition matrices $G_{i,j}$
with regular coefficients, i.e. rational functions that have no pole in $U_i\cap U_j$.
This is Proposition $18$ of GAGA \cite{MR0082175}. 
\item If we have a meromorphic connection $\tilde{\nabla}$ on $\tilde{V}$, in such trivializations, the entries of the connection matrices of $\hat{\nabla}$ are rational one forms.
Indeed, $G_{ij}$ extends to a global rational function to $\mathrm{GL}_m(\C)$, and fixing $H_{i_0}\equiv Id$, $H_i=G_{i,i_0}$ we have a global birational trivialization: $G_{i,j}=H_iH^{-1}_j$. In this trivialization, the connection is given by a global connection matrix and its coefficients must be rational one forms.
\end{enumerate}
Also  observe that, if $(\tilde{V},\tilde{\nabla})$ is a $\mathfrak{sl}_m(\C)$-connection, so can be chosen its algebraization, by \cite[Th\'eor\`eme $3$]{SerreChev58}.
\end{remark}
We will prove the following slightly more precise version of Theorem $A$.
\begin{theorem}
Let $\nabla$ be a logarithmic connection on a rank $m$ holomorphic vector bundle over $\Pu$ with $n$ poles $x_1,\ldots,x_n$, $n\geq 4$.  

Let $\rho : \pi_1(\Pu\setminus \{x_1,\ldots,x_n\},x_{0})\rightarrow \mathrm{GL}_m(\C)$ be its monodromy representation.
Suppose the germs of $\nabla$ at the poles are mild transversal models and $\rho$ is semisimple or $m=2$.
Then the following are equivalent.
\begin{enumerate}
\item \label{unthm} The conjugacy class $[\rho]$ has a finite orbit under the pure mapping class group $\mathrm{P}\M_n(\Pu)$.
\item  \label{deuxthmlog} The germ of universal isomonodromic deformation of $\nabla$ is logarithmically algebraizable.
\item \label{deuxthm} The germ of universal isomonodromic deformation of $\nabla$ is algebraizable.
\end{enumerate}
\end{theorem}
\begin{proofbis} We have seen in section $\ref{BraidMCG}$, that $[\rho]$ has finite orbit under $\mathrm{P}\M_n(\Pu)$ if and only if it has finite orbit under the pure braid group $\pi_1(F_{0,n}\Pu,x)$.~\\

 The assertion $\ref{deuxthmlog}. \Rightarrow \ref{deuxthm}.$ is trivially true.

\noindent $\ref{deuxthm}. \Rightarrow \ref{unthm}$. With the notation of Definition \ref{defalg}.
Let $D$ be the image under $\pi$ of the polar locus of $\tilde{\nabla}$. Let $D^{v}$ be the union of the irreducible components of $D$ that don't dominate $Y$.
Let $A\subset Y$ be the image by $q$ of the exceptional locus of $\pi:X\rightarrow Y\times \Pu$. It is an analytic subspace of codimension at least $1$.
Let $U:=Y\setminus (A \cup q(D^{v}))$, the conditions of Definition $\ref{defalg}$ ensure that $\psi(\Delta)$ is a neighborhood of $\psi(b_0)$  contained in $U$.

  Restricting over $U$ we have the commutative diagram below.
 $$\xymatrix{
\Delta \times \Pu \ar[r]^-{\Psi} \ar[d]_{p}&X_{\vert U} \ar[d]^{q}\ar[r]_-{\pi}^-{\simeq}&U\times \Pu \ar[d]^{\tilde{q}}\\
    \Delta \ar[r]^{\psi} & U\ar@{=}[r]&U}$$
    Let $D^h=D\setminus D^{v}$, every component $D_i$ of $D^h$ defines a finite branched covering  $q_i: D_i \rightarrow Y$, by restriction of $q$. Every $q_i$ is \'etale in the neighborhood of $\psi(b_0)$, because $\Psi$ lifts to a connection isomorphism in the neighborhood of $b_0$. 
There exists a generically finite surjective map  $u:Y_0\rightarrow Y$, \'etale in the neighborhood of  $\psi(b_0)$, such that in the pullback bundle $Y_0\times \Pu\rightarrow Y_0$, every irreducible component of $(u \times id_{\Pu})^{-1}D^h$ meets every fiber exactly once. Existence of such a map $u$ is proved in \cite[Theorem $18.2$ p $56$]{MR2030225}.

By pullback \textit{via} $u$, provided $\Delta$ is small enough, we obtain the existence of a Zariski open subset $W$ of $Y_0$ and an open embedding $\Psi_0: \Delta \times \Pu \rightarrow V \times \Pu$ with a commutative diagram
\[\xymatrix{\Delta \times \Pu \ar[r]^{\Psi_0} \ar[d]_{p}&W\times \Pu \ar[d]^{\tilde{q}_0}\\
    \Delta \ar[r]^{\psi_0}&W}\]
    and an algebraic flat meromorphic connection $\tilde{\nabla}_0$ on $W\times \Pu$ with polar locus given by disjoint sections of $W\times \Pu \rightarrow W$ such that $\Psi_0^*\tilde{\nabla}_0$ is isomorphic to $\hat{\nabla}_{\vert \Delta}$.

    By algebraicity, the polar sections take the form $v\mapsto (v,g_i(v))$ for $n$ finite maps $(g_i)_{i=1,\ldots,n}$. These functions $g_i$ are coordinate functions for a finite map $g: W\rightarrow F_{0,n}\Pu$. The map $g$ is dominant, because the restriction $\Delta\rightarrow F_{0,n}\Pu$ of the covering map $B_0\rightarrow F_{0,n}\Pu$ factors through it.
   By \cite[Lemma 4.19]{MR1841091}, we can deduce that the image of $g_* :\pi_1(V,\psi_0(b_0))\rightarrow \pi_1(F_{0,n}\Pu,x)$ is a finite index subgroup of the pure braid group $\pi_1(F_{0,n}\Pu,x)$. By Corollary $\ref{finiteindex}$, this allows to conclude that $[\rho]$ has finite orbit under $\pi_1(F_{0,n}\Pu,x)$.    ~\\

 \noindent $\ref{unthm}. \Rightarrow \ref{deuxthmlog}$. We use Theorem $\ref{topthm}$ and Theorem $\ref{topthm rang2}$. With the notation of these latter theorems, let $Y$ be a a projective manifold in which $U$ embeds as a Zariski open subset.
    The total space $E$ of $f^*\mathcal{T}_n$ embeds naturally in $Y\times \Pu$ as the complement of an algebraic subspace $A$, let $D$ be the codimension $1$ part of $A$. Let $\pi:X\rightarrow Y\times \Pu$ be an embedded desingularization of $D$, which is an isomorphism on its image in the neighborhood of $\{b\}\times \Pu$ , set $D^0:=\pi^{-1}(D)$.

    We will realize the representation $\pi^*\hat{\rho}$ as the monodromy representation of a well chosen logarithmic flat connection with polar locus in $D^0$. Let $D_1,\ldots, D_n$ be the components of $D^0$ that intersect the fiber $q^{-1}(b)$ of the ruling $q:X\rightarrow Y$, in order to have
     $\pi(D_i\cap q^{-1}(b))=(b,x_i)\in Y\times \Pu$. In the product $U\times \Pu$ the hypersurface $\pi(D_i)$ is simply the graph of the coordinate function $f_i$ of $f$.

 By Theorem $\ref{RHthm}$, we can define a realizing flat logarithmic connection $\tilde{\nabla}$ over $X$ for $\pi^*\hat{\rho}$ by assignment of a compatible mild transversal type to each component of $D^0$.
 We now describe our choice. For $D_i$, take the mild transversal type defined by the germ of $\nabla$ at its pole $x_i$; for components that do not intersect $q^{-1}(b)$, choose any mild transversal model compatible with $\pi^*\hat{\rho}$.

 By the uniqueness part of Theorem $\ref{RHthm}$, the pullback of $\nabla$ by the composition $X\rightarrow Y\times\Pu \rightarrow \Pu$ is isomorphic to the restriction of $\tilde{\nabla}$ to $q^{-1}(b)$.
 Let $\Delta_U$ be a (euclidean) simply connected neighborhood of $b$ in $U$ on which $f$ is an analytic isomorphism on its image $\Delta_1\subset F_{0,n}\Pu$. If $\Delta_U$ is chosen small enough, we have also an isomorphism $(\Delta,b_0)\rightarrow(\Delta_1,x)$ induced by $B_0\rightarrow F_{0,n}\Pu$ and 
we obtain, by composition, an isomorphism $\psi : (\Delta,b_0)\rightarrow (\Delta_U,b)$, for a euclidean neighborhood $\Delta_U$ of $b$ in $U$. The inverse of $\psi$ is a local lift of $f$, $\tilde{f}:(\Delta_U,b)\rightarrow (B_0,b_0)$.

We have a commutative diagram.  $$\xymatrix{
\Delta \times \Pu \ar[r]_-{\psi\times id} \ar[d]_{p}  \ar@/^1.2pc/[rr]^{\Psi:=}&\Delta_U\times \Pu \ar[d]^{\tilde{q}}&X_{\vert \Delta_U} \ar[d]^{q}\ar[l]_-{\simeq}^-{\pi}\\
    \Delta \ar[r]^{\psi} &\Delta_U& \Delta_U\ar@{=}[l]}$$
Application of Proposition $\ref{propuniv}$ gives an isomorphism between $\pi_* \tilde{\nabla}_{\vert X_{\vert \Delta_U}}$ and $(\psi\times id)_*\hat{\nabla}_{\vert \Delta\times \Pu}$. This implies that we have an isomorphism \[(\hat{V},\hat{\nabla})_{\vert \Delta\times \Pu} \simeq \Psi^*(\tilde{V},\tilde{\nabla}),\] for the map $\Psi$ defined in the diagram.
We see that every condition of Definition~$\ref{defalg}$ for logarithmic algebrization is fulfilled.
\qed \end{proofbis}

\begin{remark}
For  $\ref{deuxthm}. \Rightarrow \ref{unthm}.$ we \emph{did not use} the semisimplicity assumption.
\end{remark}

\begin{remark}
In view of the last statement of Theorem $\ref{topthm}$, in the proof of $\ref{unthm}. \Rightarrow \ref{deuxthmlog}.~$,
if the initial connection $\nabla$ splits as $\nabla=\oplus_{i\in I}\nabla_i$,
 the logarithmic  algebraization $\tilde{\nabla}$ we construct has the same splitting type $\tilde{\nabla}=\oplus_{i\in I}\tilde{\nabla}_i$.\end{remark}
 \begin{remark}
 In \cite{MR2289083} Inaba, Iwasaki and Saito show that the  moduli spaces of $\alpha$-stable rank $m$ $\lambda$-parabolic connections over $\Pu$ corresponding to different positions of the $n$ poles arrange in an algebraic family
 $\pi :M^{\alpha}(\lambda)\rightarrow T_n'$, where $T_n'$ is a finite \'etale cover of $F_{0,n}\Pu$.
  Then, they prove that the Riemann-Hilbert correspondance allows to define a non singular holomorphic foliation $\F_{\alpha}$ of dimension $n$ on $M^{\alpha}(\lambda)$, transverse to the fibration $\pi$, corresponding to isomonodromy.

  It seems likely that our Theorem \ref{maintheorem} would allow to show that any  $\alpha$-stable rank $m$ $\lambda$-parabolic connection over $\Pu$ whose monodromy representation yields a finite braid group orbits on the suitable character variety  corresponds to an algebraic leaf of the foliation $\F_{\alpha}$.
 \end{remark}
 Relaxing the conclusion, we can remove the mildness hypothesis in Theorem~\ref{maintheorem}.
 \begin{COR:1}
 Let $\nabla$ be a logarithmic connection on a rank $m$ holomorphic vector bundle over $\Pu$ with $n$ poles $x_1,\ldots,x_n$, $n\geq 4$.  

Let $\rho : \pi_1(\Pu\setminus \{x_1,\ldots,x_n\},x_{0})\rightarrow \mathrm{GL}_m(\C)$ be its monodromy representation.
Suppose  $\rho$ is semisimple or $m=2$.
Then the following are equivalent.
\begin{enumerate}
\item  The conjugacy class $[\rho]$ has a finite orbit under the pure mapping class group $\mathrm{P}\M_n(\Pu)$.
\item   Up to a birational gauge transformation of $\nabla$, the germ of universal isomonodromic deformation of $\nabla$ is algebraizable.
\end{enumerate}
\begin{proofbis}
Any logarithmic flat connection is birationally gauge equivalent to a flat logarithmic connection with mild transversal types (\textit{e.g.} elements of $Eul^0$). 
\qed \end{proofbis}
 \end{COR:1} 
\subsection{Garnier systems}\label{secGarnier}
We now want to explain an application of our main theorem for rank two connections. It will relate finite braid group orbits and algebraic solutions of a family of isomonodromy equations, namely Garnier systems.

\subsubsection{Garnier divisor}

Let $\nabla$ be a logarithmic trace free connection on $\mathcal{O}^{\oplus 2}_{\Pu}$, with $n=N+3$ non-apparent distinct poles $t_1^0,\cdots,t_{N}^0,0,1,\infty$. 
Non-apparent means the monodromy of $\nabla$ around any of these poles is not scalar.
Trace free means that the connection matrix of $\nabla$ is trace free.
In the sequel, we will assume that our base point $x$ for $F_{0,n}\Pu$ is $x=(t_1^0,\ldots,t_N^0,0,1,\infty)$.
Let $(\hat{V},\hat{\nabla})$ be the universal isomonodromic deformation of $\nabla$.

Recall that we denote $r: (B_0,b_0)\rightarrow (F_{0,n}\Pu,x)$ the universal cover of $F_{0,n}\Pu$ and $F_{3,N}\Pu\subset F_{0,n}\Pu$ is the subset of normalized configurations $F_{3,N}\Pu=\{(t_1,\ldots,t_N,0,1,\infty)\in F_{0,n}\Pu\}$.
Define $C_0:=r^{-1}(F_{3,N}\Pu)$, recalling  $F_{0,n}\Pu\simeq F_{3,N}\Pu\times \mathrm{PSL}_2(\C)$, we see that $r_{\vert C_0}: C_0\rightarrow F_{3,N}\Pu$ is a universal covering.

 We will take special interest in the restriction $(\hat{V}_{\vert C_0 \times \Pu},\hat{\nabla}_{\vert C_0\times \Pu})$.
Let $D_i=(r\times id_{\Pu})^{-1}(\{(y,z)\in F_{3,N}\Pu \times \Pu \vert z=y_i\})$.
Let $R_i\in \mathrm{H}^0(\mathrm{End}(\hat{V}_{\vert D_i}),D_i)$ be the residue of $\hat{\nabla}_{\vert C_0 \times \Pu}$ on $D_i$, it is defined locally by the residue of the connection matrix. The conjugacy class of $R_i(p)$ is independent of $p\in D_i$. Let $\theta_i$ be such that the spectrum of $R_i$ is $\{\theta_i/2,-\theta_i/2\}.$

As the pole $D_i$ is non-apparent, we have a unique holomorphic section $v_i^{\theta_i}$ of $\mathrm{P}(\hat{V}_{\vert D_i})$ that gives the direction of the eigenspace of $R_i$ associated to $\theta_i/2$.
Define $K:=\{y\in C_0\vert  ~\hat{V}_{\vert \{y\} \times \Pu} \mbox{ is not trivializable}\}$. It is an analytic hypersurface of $C_0$ ($\Theta$-divisor) \cite{MR0728431}.

For any $y\in C_0\setminus K$, we have a unique holomorphic  section $\sigma^y$ of $\mathrm{P}(\hat{V}_{\vert \{y\} \times \Pu})$ that coincides with $v_n^{\theta_n}$ over ${(y,\infty)}$ and is given by a constant map $\Pu\rightarrow \Pu$ in any trivialization $\mathrm{P}(\hat{V}_{\vert \{y\} \times \Pu})\simeq \Pu\times \Pu$.

Here, one needs to recall a piece of foliation terminology.
The foliation of the total space of $\hat{V}_{\vert \{y\} \times \Pu}$  given by the horizontal sections of $\hat{\nabla}_{\vert \{y\} \times \Pu}$ is invariant under the natural $\C^*$ action.
The quotient foliation on $\mathrm{P}(\hat{V}_{\vert \{y\} \times \Pu})$ is called a Riccati foliation and denoted  $\mathrm{P}(\hat{\nabla}_{\vert \{y\} \times \Pu} )$.
For a holomorphic foliation $\F$ on a smooth complex  surface $S$, for any embedded curve $\mathcal{C}\subset S$ not invariant by $\F$, the tangency divisor of $\mathcal{C}$ with $\F$ is the divisor $\sum_{p\in \mathcal{C}} n_p p$ on $\mathcal{C}$, with $n_p:=dim_{\C}\mathcal O_{S,p}/I_p$ where $I_p$ is the ideal $<f,v(f)>$, for $f$ a local equation of $\mathcal{C}$  at $p$, and $v$ a local vector field tangent to $\F$, with isolated singularities.

Notice that Condition $\ref{noninv}$ below is guaranteed by irreducibility of $\nabla$ or, using Camacho-Sad formula, by  $\left [ \forall (\varepsilon_i)\in \{+1,-1\}^n,\sum \varepsilon_i \theta_i\neq 0\right]$.
\begin{cond}\label{noninv}
The image of $\sigma^y$ is not an invariant curve for the foliation $\mathrm{P}(\hat{\nabla}_{\vert \{y\} \times \Pu} )$ on $\mathrm{P}(\hat{V}_{\vert \{y\} \times \Pu})$.
\end{cond}
Under this condition, we can set the following.
\begin{definition} We define a divisor $G_y=G_y(\nabla,\theta_n)$ on $\{y\}\times \Pu$ by requiring $(G_y+\{(y,\infty)\})$ to be the projection of the tangency divisor between the Riccati foliation $P(\hat{\nabla}_{\vert \{y\} \times \Pu} )$ and the curve $\sigma^y(\{y\}\times \Pu)$.
\end{definition} 

\begin{prop/def}
If Condition $\ref{noninv}$ is fulfilled for general $y\in C_0\setminus K$,
there exists a divisor $G=G(\nabla,\theta_n)$ on $(C_0\setminus K) \times \Pu$ such that, for every $y\in C_0\setminus K$ satisfying Condition $\ref{noninv}$,  $G_y=G_{\vert \{y\} \times \Pu}$.

\textnormal{This divisor $G$ is called \emph{the Garnier divisor} of $\nabla$ with respect to $\theta_n$.}
\end{prop/def}
\begin{proofbis} By \cite{MR0728431}, for any point $y\in C_0\setminus K$, there exists a neighborhood $\Delta_0$ of $y$ in $C_0 \setminus K$ and a bundle isomorphism
 $\phi :\mathcal{O}^{\oplus 2}_{\Delta_0\times \Pu}\stackrel{\simeq}{\longrightarrow} \hat{V}_{\vert \Delta_0\times \Pu}$.
Hence, $\phi^*\hat{\nabla}_{\vert \Delta_0\times \Pu}$ has the form $Z\mapsto dZ-(Qdz+\Omega)Z$, for matrices $Q$ and $\Omega$ as below.
\begin{equation}\label{systemdecomp} \left \lbrace\begin{array}{l}Q=\sum\limits_{i=1}^N \frac{A_{i}(t)}{z-t_i}+\frac{A_{N+1}(t)}{z}+\frac{A_{N+2}(t)}{z-1};\\
\Omega=\sum\limits_{i=1}^N\Omega_i(x,t) dt_i;
\end{array} \right .
\end{equation}
where $(t_i)$ are the obvious local coordinates on $C_0$ induced by $r$.

 Define $A_{\infty}=-\sum_{i=1}^{N+2} A_i$. For $i=1,\ldots,n$, the matrix function $t \mapsto A_i(t)$ gives the residue of $\phi^*\hat{\nabla}_{\vert \Delta_0\times \Pu}$ on $D_i$.
If $\phi$ is well chosen, for any $t$, the eigenline of $A_{\infty}(t)$ associated with $\theta_{n}/2$ is spanned by $(1,0)$, that is $$A_{\infty}(t)=\left [\begin{smallmatrix} \theta_{n}/2 & *\\ 0&-\theta_n/2\end{smallmatrix} \right ].$$

With such a choice, if we write $Q_{2,1}(z,t)=\frac{c(z,t)}{z(z-1)\prod_{i=1}^N(z-t_i)}$, a simple computation shows that $G_y$ is the degree $N$ divisor on $\{y\}\times \Pu$ that coincides with the zero divisor of $z\mapsto c(z,t(y))$ on $(\{y\}\times \Pu)\setminus \{(y,\infty)\}$.
For this reason, $G$ can be defined by an equation $c_{hom}(z,\xi,t)=0$ on $\Delta_0\times \Pu$, for $c_{hom}$ a homogeneous polynomial of total degree $N$ in the variables $z,\xi$ and with coefficients given by holomorphic functions on $\Delta_0$.
\qed \end{proofbis}
At this point it is natural to ask if the Garnier divisor is reduced in general. 
We may also ask if it can have contact (or even share a component) with the polar hypersurface $D=\cup_i D_i$.
To the author's knowledge, there is no information on these questions in the literature.

\begin{remark}
In the local computation above, we see that Condition $\ref{noninv}$ for $y$ is fulfilled if and only if $z\mapsto c(z,t(y))$ is a non zero polynomial map. Consequently, it is tantamount to prove that this condition is fulfilled for a given $y$ or for general $y\in C_0\setminus K$.
\end{remark}

\begin{definition}
We will say that the Garnier divisor of $\nabla$ with respect to $\theta_n$ is well defined if Condition $\ref{noninv}$ is fulfilled by some $y\in C_0\setminus K$.
\end{definition}

\subsubsection{Algebraicity of the Garnier divisor}
Let $\nabla$ be a logarithmic trace free connection on $\mathcal{O}^{\oplus 2}_{\Pu}$, with $n=N+3$ non-apparent distinct poles $t_1^0,\ldots,t_{N}^0,0,1,\infty$. 
\begin{theorem}\label{AlgGD}
If the Garnier divisor $G(\nabla,\theta_n)$ is well defined and the germ of universal isomonodromic deformation of $\nabla$ is algebraizable, then
the projection by $(C_0\setminus K) \times \Pu\rightarrow F_{0,N}\Pu\times \Pu$ of the support of the Garnier divisor of $\nabla$ is a union of graphs of algebraic functions from $F_{0,N} \Pu$ to $\Pu$ \textit{i.e.} its closure in $\left(\Pu\right)^N\times \Pu$ is an algebraic hypersurface.
\end{theorem}

\begin{proofbis}
By assumption, $\nabla$ meets the conditions of Definition $\ref{defalg}$ and we use the corresponding notation.
Let $U$ be a Zariski neighborhood of $\psi(b_0)$ in $Y$ over which $X\rightarrow Y$ is isomorphic, \textit{via} $\pi$, to the trivial fibration $U\times \Pu \rightarrow U$.

Consider the component $D_{i}$ of the polar locus of $\pi_*\tilde{\nabla}$ that passes through $(\psi(b_0),x_i)\in U\times \Pu$.
Up to a a finite covering of $U$, we may suppose every $D_i$ is the image of a regular local section $f_i:U\rightarrow U\times \Pu$ of $q$. The residue $R$ of $\pi_*\tilde{\nabla}$ on $D_{n}=f_n(U)$ is a regular section of $\mathrm{End}(\tilde{V})_{\vert f_n(U)}\rightarrow f_n(U)$. The eigendirection of $R$ corresponding to  $\theta_n$ gives a regular section $\sigma_0$ of $\mathbb{P}(\tilde{V})_{\vert f_n(U)}$ and we may apply Lemma $\ref{sectionreg}$ below, with $s=f_n$ and $x=\psi(b_0)$. In this way, we obtain a regular extension $\sigma$ of $\sigma_0$ defined on $U_0\times \Pu$, with $U_0$ a Zariski neighborhood of $x$ in $U$ such that, for every $y\in U_0$,  $\mathbb{P}(\tilde{V})_{\vert \{y\}\times \Pu}$ is trivializable.

Let $U_1\subset U_0\times \Pu$ be a Zariski open set such that $\mathbb{P}(\tilde{V})_{\vert U_1}$ is algebraically trivializable and fix a trivialization.
Let $Z\mapsto dZ -\left [\begin{smallmatrix} \beta/2 &\alpha\\ -\gamma&-\beta/2\end{smallmatrix} \right ]\cdot Z$ be the local form of $\nabla$ and  let $(y,z)\mapsto h(y,z)$ be the local (regular) form of $\sigma$.
Consider the splitting $\Omega^1_{U_0\times \Pu }=\mathcal{D}_{U_0}\oplus \mathcal {D}_{\Pu}$ induced by the product structure and denote $p_z: \Omega^1_{U_0\times \Pu }\rightarrow \mathcal {D}_{\Pu}$ the projection to the rank $1$ factor.
The intersection of $U_1$ with the projection of the tangency locus between $\mathbb{P}(\tilde{\nabla})_{\vert {y}\times \Pu}$ and $\sigma_{\vert {y}\times \Pu}$ is given by the solutions $(y,z)\in U_1 \subset U_0\times \Pu$ of
$p_z(-dh+\alpha + h \beta + h^2 \gamma)=0$. We see that this locus is obtained by restriction to $\{y\}\times \Pu$ of an algebraic divisor $\tilde{G}$ on $U_0\times \Pu$.

A \emph{fundamental observation} is that $\tilde{G}$ pulls back to $G_{\vert (C_0\cap \Delta) \times \Pu}$ \textit{via} the map $\psi_{\vert (C_0\cap \Delta)}\times id_{\Pu}$.
Let $f: U_0 \rightarrow F_{0,n}\Pu$ be the map with coordinate functions $(f_i)_i$.
By Chevalley's Theorem, the image $G_0$ of $supp(\tilde{G})\subset U\times \Pu$ in $F_{0,n}\Pu\times \Pu$ by the map $f\times id_{\Pu}$ is a constructible subset of $(\Pu)^n\times\Pu$.
For this reason, its euclidean closure $\overline{G_0}$ in $(\Pu)^n\times\Pu$ is a hypersurface. As every component of $\tilde{G}$ dominates $U$, every component of $\overline{G_0}$ dominates the first factor of
$(\Pu)^n\times\Pu$. 

We now consider the image $G_1$ of the support of the Garnier divisor $G=G(\nabla,\theta_n)$ by \[r_{\vert C_0\setminus K} \times id_{\Pu} : (C_0\setminus K) \times \Pu\rightarrow F_{0,n}\Pu\times \Pu.\]
By restriction of $r_{\vert C_0\setminus K} \times id_{\Pu}$, for each component $\kappa$ of $G$ we have an analytic map $ L: \kappa  \rightarrow F_{3,N}\Pu\times \Pu \subset F_{0,n}\Pu\times \Pu$
let $w$ be a point of $\kappa$ lying over $b_0$ in  $(C_0\setminus K) \times \Pu$.
Let $g$ be a local (algebraic) equation  for $G_0$ in a Zariski  neighborhood $W$ of $(r_{\vert C_0\setminus K} \times id_{\Pu}) (w)$ in $F_{0,n}\Pu\times \Pu$.
In a euclidean neighborhood of $w$ in $\kappa$, we have $g\circ L \equiv 0$, by our fundamental observation. By analytic continuation, this holds on the whole of $L^{-1}(W)$ and $G_1$ is contained in $\overline{G_0}$.

Again by the fundamental observation, every component of $\overline{G_0}\cap(F_{3,N}\Pu\times \Pu)$ is met by $G_1$ in a euclidean open subspace.
By \'etaleness of  $r$, we see that $G_1$ is a closed analytic subspace of $\tilde{G}_0=G_0\cap r(C_0\setminus K)$, with equal dimension. As it meets every of the components of $\tilde{G}_0$ in an open set, we may conclude $G_1=\tilde{G}_0$ and  $\overline{G_1}$= $\overline{G_0}$, which completes the proof.



\qed \end{proofbis}

\begin{lemma}\label{sectionreg}
Let $V\rightarrow U\times \Pu$ be a rank $2$ algebraic vector bundle, $U$ being a quasi-projective manifold.
Denote $q : U\times \Pu \rightarrow U$ the projection to the first factor.

Let $s: U\rightarrow U\times \Pu$ be a regular section of $q$ and $\sigma_0$ be a regular section of $\mathbb{P}(V)_{\vert s(U)}$.
Suppose that there exists $x \in U$ such that $\mathbb{P}(V)$ is trivializable in restriction to $q^{-1}(x)$.

Then 
\begin{enumerate}
\item \label{Maruyama} There exists a Zariski neighborhood $U_0$  of $x$ in $U$ such that, for any $y\in U_0$, the restriction of $\mathbb{P}(V)$ to $q^{-1}(y)$ is trivializable.
\item \label{regular}This allows to define a set theoretic section $\sigma$ of $\mathbb{P}(V)_{\vert U_0}$, by requiring that, for every $y\in U$,  $\sigma_{\vert q^{-1}(y)}$ is constant in any trivialization  of $\mathbb{P}(V)_{\vert q^{-1}(y)}$ and coincides with $\sigma_0$ over $s(U_0)$. This section is actually regular.
\end{enumerate}
\end{lemma}
\begin{proofbis}
The assertion $(\ref{Maruyama})$ is direct application of the openness result of Maruyama \cite{MR0429899}, because projectively trivial vector bundles on $\Pu$ are those of type $(0,\ldots,0)$.
In $(\ref{regular})$, the nontrivial part is the regularity of $\sigma$. Without loss of generality, we assume that $q^{-1}(U_0)$ does not contain singular fibers of $q$. In the proof, we denote $P_y:=\mathbb{P}(V)_{\vert q^{-1}(y)}$, for every $y \in U_0$.
We have a bundle isomorphism $P_y\simeq \mathbb{P}^{1}\times \Pu$.
Let $S_y$ be the set of curves in $P_y$ that map to graphs of constant sections in this isomorphism. Let $S=\cup_{y\in U_0} S_y$.
We will show that the elements of $S$ are the level sets of a morphism of quasi-projective varieties $W\rightarrow Q$, where $W$ is the total space of $\mathbb{P}(V)_{\vert U_0}$.

For this purpose, we will use a piece of deformation theory for morphisms from $\Pu$ to algebraic varieties, see \cite[chap. $2$ and $3$]{MR1841091}.

 Let $C\in S$, and choose an embedding $f :\Pu \rightarrow W$ with image $C$. the result \cite[Prop. $4.1$]{MR0728431} shows that we can trivialize holomorphically $V$ in an euclidean neighborhood of any $q^{-1}(y),$ $y\in U_0$. 
 
For this reason,  the map $f$ satisfies $$f^*TW=T\Pu\oplus \mathcal O_{\Pu}^{\oplus \dim W-1}.$$
This implies that $[f]$ is a smooth point of $\mathrm{Mor}(\Pu,W)$ and that the irreducible component $M$ of $\mathrm{Mor}(\Pu,W)$ passing through $[f]$ has dimension $\dim W+2$. 

Let $\overbull{M}$ be the Zariski dense subset of $M$ consisting in morphisms $g:~\Pu\rightarrow~W$ with $$g^*TW=T\Pu\oplus \mathcal O_{\Pu}^{\oplus \dim W-1}.$$
The above Zariski density is again obtained by the openness result of Maruyama.

Let $M_0$ be the closed subset of $M$ of elements $g\in M$ such that $ pr \circ g=id_{\Pu}$ where $pr$ is the composition of natural maps $pr : W\stackrel{r}{\rightarrow} U_0\times \Pu\rightarrow \Pu$. We will also use the notation $\overbull{M}_0:=M_0\cap \overbull{M}$. If $f$ is well chosen, $f \in M_0$. We always assume $f\in M_0$ in the sequel.

Firstly, consider the compositions $r\circ g$, $[g]\in M$  where $r$ is the bundle projection $W\stackrel{r}{\rightarrow} U_0\times \Pu$.
As above, the morphism $[r\circ f]$ is  a smooth point of $\mathrm{Mor}(\Pu,U_0\times \Pu)$, and the component through $[r\circ f]$ is isomorphic to $m=\mathrm{PSL}_2(\C) \times U_0$, the evaluation map being \[\begin{array}{ccc} \mathrm{PSL}_2(\C) \times U_0\times \Pu& \longrightarrow &U_0\times \Pu \\
					(A,y,z)&\longmapsto &(y, A\cdot z).\end{array}\]
By universality, we have a morphism $\psi: M\rightarrow m$  such that, for any $g\in M$, $\psi([g])=[r\circ g]$. In particular any of the images $g(\Pu)$, $[g]\in M_0$ is contained in $P_y$, for a certain $y \in U_0$ and projects isomorphically onto $q^{-1}(y)$ by $r$. 

Secondly, every element $[g]\in \overbull{M}_0$ satisfies $g(\Pu)\in S$. Indeed, by the previous point $g(\Pu)\subset P_y$, for a certain $y\in U_0$, and the normal bundle of $g(\Pu)$ in $P_y$ is trivial (as any subbundle of a trivial bundle over $\Pu$). Under $P_y\simeq \mathbb{P}^{1}\times \Pu$, $g(\Pu)$ is the graph of a morphism $\Pu \rightarrow \mathbb P^{1}$. The triviality of its normal bundle is equivalent to this morphism being constant and, in turn, to $g(\Pu)\in S_y \subset S$.

Thirdly, the restriction of the evaluation map $\overbull{M}_0 \times \Pu\rightarrow W$ is dominant: this can be obtained by injectivity of this restriction and the estimate $\dim M_0\geq \dim M-3$, thanks  to Brouwer's invariance of domain. The estimate comes from the observation that, for an element $g\in M$, the condition $ pr \circ g=id_{\Pu}$ is tantamount to the conjunction of the three identities $ pr \circ g(0)=0$,  $ pr \circ g(1)=1$, $ pr \circ g(\infty)=\infty$. This implies that the composition $\overbull{M}_0 \times \Pu\rightarrow W\rightarrow U_0$ is dominant.

Then, for any element $[g]\in \overbull{M}_0$, with  $r\circ g(\Pu)=q^{-1}(y)$, $[g]$ is a smooth point of $\mathrm{Mor}(\Pu,W)$, but we can also interpret it as a smooth point of $\mathrm{Mor}(\Pu,P_y)$, and by 
universality, every element of $\mathrm{Mor}(\Pu,P_y)$, contained in the component of $[g]$ can be seen as a point of $\overbull{M}$. In particular, any element of $S_y$ has the form $\tilde{g}(\Pu)$, for some $[\tilde{g}]\in \overbull{M}_0$. 

Finally, we see that, upon restriction of the initial Zariski neighborhood $U_0$ of $x$, we have a bijective morphism $\overbull{M}_0 \times \Pu\rightarrow W$, by Zariski's Main Theorem, this is an isomorphism and we get a fibration $W\rightarrow \overbull{M}_0$ with fibers equal to the elements of $S$.

The graph of the section $\sigma$ is the saturation of $\sigma_0$ by fibers of this fibration, it is thus a closed Zariski subset of $W$; this proves that $\sigma$ is a regular section of $\mathbb{P}(V)_{\vert U_0\times \Pu}\rightarrow U_0\times \Pu$. The reduction of $U_0$ in the proof is harmless, because we can apply the above proof with any initial base point $x$.
\qed \end{proofbis}

\subsubsection{Conclusion}
We now explain the relation between Garnier divisor, Garnier system and isomonodromy.
Let $\nabla_{t^0}$ be a trace free logarithmic connection on  $\mathcal{O}^{\oplus 2}_{\Pu}$, with $n=N+3$ non-apparent poles
given by the coordinates of our base point $x=t^0=(t_1^0,\ldots,t_N^0,0,1,\infty)$.
Let $T$ be a simply connected neighborhood of $t^0$ in $F_{3,N}\Pu$ and consider the germ of universal isomonodromic deformation of $\nabla_{t^0}$ restricted to $T$---formally we are identifying $T$ with a lift in the universal cover of $F_{3,N}\Pu$.
For $T$ small enough, this deformation is given by  \[\nabla_t :Z\mapsto dZ-Q(z,t)dz\cdot Z\]

Let $\nabla$ be the corresponding flat logarithm connection on $\mathcal{O}^{\oplus 2}_{T\times\Pu}$.
Up to a holomorphic gauge transformation, we may suppose the eigendirection $\delta_{\theta_n}$ with eigenvalue $\theta_n/2$ of the residue of $\nabla$ over $z=\infty$ is given by the constant section
$(1,0)$ of $(\mathcal{O}^{\oplus 2}_{T\times\Pu})_{\vert T\times \{\infty\}}$.
 
 Let $a(z,t),b(z,t),c(z,t)$ give the entries of $Q$ as follows, $Q=\left [\begin{smallmatrix}b/2&a\\-c&-b/2\end{smallmatrix}\right]$.
 
We suppose that the following condition is fulfilled.
\begin{cond}\label{condgf}
 the Garnier divisor $G=G(\nabla_{t^0},\theta_n)$ is well defined, shares no component with the polar locus of $\nabla$ and $G_{\vert \{t^0\}\times \Pu}$ is reduced.
\end{cond} 
Under this condition, up to shrinking $T$, we have local equations of the form
$z=\lambda_i(t)$, $i=1,\ldots,N$ 
for the components of $G$ in $T\times \Pu$, with holomorphic $\lambda_i: T\rightarrow \Pu$. In this way, $c$ writes
$c(z,t)=c_0\psi/\phi$ where $c_0\in \C^*$, $\psi=\prod_{i=1}^N (z-\lambda_i(t))$ and $\phi=z(z-1)\prod_{i=1}^N  (z-t_i)$.
Also, by logarithmicity $a=q_a(z,t)/\phi(z,t)$ and $b=q_b(z,t)/\phi(z,t)$ with $q_a,q_b$ polynomials in $z$ of degree $\leq N+1$.

Using a holomorphic gauge transformation given by a matrix $\left [\begin{smallmatrix}\lambda(t)&0\\0&1/\lambda(t)\end{smallmatrix}\right]$, we obtain
$c_0=1$. In addition, using a second one, given by a matrix of the form$\left [\begin{smallmatrix}1&g(t)\\0&1\end{smallmatrix}\right]$, we obtain
that the $1$-form $b(z)dz$ expands with no constant term at $z=\infty$. Namely, for $\xi=1/z$, in the neighborhood of $\xi=0$,   $b(z)dz=(\frac{-\theta_n}{\xi}+o(1))d\xi$.
\begin{definition}
With such $c$ and $b$, we will say that $(\nabla_t)$ is in \emph{normalized form with respect to $\theta_n$.}
\end{definition}
Let $H(z,t):=\left [\begin{smallmatrix}0&-1/\sqrt{c}\\ \sqrt{c}&\tau \sqrt{c}\end{smallmatrix}\right]$, where $\tau=\frac{1}{2c}(\frac{\partial_zc}{c}+b)$.
The transformation $\mathcal T: Z\mapsto H(z,t)\cdot Z$ is well defined only on a  ramified covering of $\Pu\times T$, but the induced transformation on the level of projectivized connections is well defined, in particular the transformed connection $\mathcal T_*\nabla$ descends
on $\Pu\times T$ and is flat. Computation shows that the family $\mathcal T_*\nabla_t$ is given by $\tilde{Z}\mapsto \left [\begin{smallmatrix} 0&1\\ p(z,t)&0\end{smallmatrix}\right] dz\cdot \tilde{Z}$ where
 \begin{equation}\label{p1} p(z,t)=\frac{3}{4}\left(\frac{\partial_z c}{c}\right)^2-\frac{\partial^2_z c}{2c}-\frac{\partial_z b}{2}+\frac{b}{2}\frac{\partial_z c}{c}+\frac{b^2}{4}-ac.
 \end{equation}
 
 With the specific forms for $a,b,c$  given above we can see that the equation $\partial^2_zw=p(\cdot,t)w$ is Fuchsian for all $t\in T$. 
 Indeed, every term in equation $(\ref{p1})$ can be written in the form $\frac{q}{\phi^2\psi^2}$ for a polynomial $q$ in $z$ of degree $\leq 4N+2$.
 These are elementary computations, \textit{e.g.} for $b\frac{\partial_z c}{c}=b\frac{\partial_z \psi}{\psi}-b\frac{\partial_z \phi}{\phi}$, each term satisfies the required property:
\[ \begin{array}{l} b\frac{\partial_z \psi}{\psi}=\frac{q_b\psi\phi \partial_z \psi}{\phi^2\psi^2}, \deg_z q_b\psi\phi \partial_z \psi\leq (N+1)+N+(N+2)+(N-1)=4N+2,\\
 b\frac{\partial_z \phi}{\phi}=\frac{q_b \psi^2\partial_z \phi}{\phi^2\psi^2}, \deg_z q_b\psi^2 \partial_z \phi\leq (N+1)+2N+(N+1)=4N+2.\\
\end{array}\]
Defining 
\begin{equation}\label{defphipsi}
\left \lbrace \begin{array}{lc}\psi_{\lambda}(z):=\prod_{i=1}^N (z-\lambda_i)\in\C[z],&~~ \lambda\in \C^N,\\
 \phi_{t}(z):=z(z-1)\prod_{i=1}^N  (z-t_i)\in\C[z],&~~  t\in \C^N, \\
 \end{array}\right.
\end{equation}
we have the following.
\begin{lemma}\label{linalg}
Let   $0$,$1$,$(t_i)_{i=1,\ldots,N}, (\lambda_i)_{i=1,\ldots,N}$ be $2N+2$ distinct points in $\C$. The following elements  give a basis for the vector space $\frac{1}{\phi_t^2\psi_{\lambda}^2}\cdot {\C[z]}_{4N+2}$.
\[\begin{array}{c}
\begin{array}{lcr}\frac{1}{z(z-1)}~;&\left(\frac{1}{(z-t_i)^2}\right)_{i=1,\ldots, N+2} ;&\left(\frac{1}{(z-\lambda_i)^2}\right)_{i=1,\ldots, N} ; \end{array}\\
\\
\begin{array}{lr} \left(\frac{1}{z(z-1)(z-\lambda_i)}\right)_{i=1,\ldots, N} ;&\left(\frac{1}{z(z-1)(z-t_i)}\right)_{i=1,\ldots, N} .\end{array}
 \end{array}\]
\begin{proofbis} Linear algebra. \qed \end{proofbis}
\end{lemma}
By Lemma $\ref{linalg}$, we can rewrite $p$ in the following form, with uniquely determined coefficients $a_i,L_i,\nu_i$.
 \begin{eqnarray}\label{eqp}
\begin{array}{l}
p(z,t)=\frac{a_{N+1}}{z^2}+\frac{a_{N+2}}{(z-1)^2}+\frac{a_{N+3}}{z(z-1)}\\+\sum\limits_{i=1}^N \left [\frac{a_{i}}{(z-t_i)^2}+\frac{t_i(t_i-1)L_i}{z(z-1)(z-t_i)}+\frac{3}{4(z-\lambda_i)^2}-\frac{\lambda_i(\lambda_i-1)\nu_i}{z(z-1)(z-\lambda_i)}\right].
 \end{array}
 \end{eqnarray}
 The coefficients $a_i$, $i=1,\ldots,N+2$ can be determined from the local exponents $\theta_i$ by residue computations from equation $(\ref{p1})$,
 \[a_i=({\theta_i}^2-1)/4,i=1,\ldots, N+2.\] Then, computing  $\lim_{z\rightarrow \infty} z^2p$ in two manners, using expressions  $(\ref{p1})$ and $(\ref{eqp})$, we find \[a_{N+3}=\frac{1-N}{2}-\left(\sum_{i=1}^{N+2}\frac{\theta_i^2}{4}\right)+\frac{\theta_n}{2}\left(\frac{\theta_n}{2}-1\right).\] Notice that these expressions allow to recover $\theta_n$ from $(a_i)_{i=1,\ldots,n}$, up to the involution $\theta_n\mapsto -\theta_n+2$. For abstract $(a_i)_{i=1,\ldots,n}$, we denote $\theta_{n,1}(a_i),\theta_{n,2}(a_i)$ the two (possibly equal) corresponding $\theta_n$.

 ~ \\
 Conversely, starting from $p_0$ in the form  $(\ref{eqp})$ with $2N+2$ distinct poles  $(t_i), (\lambda_i)$, for a given choice of $\theta_n\in\{\theta_{n,1}(a_i),\theta_{n,2}(a_i)\}$, we can recover the unique corresponding normalized family $(\nabla_t)$. Indeed, $c=\psi/\phi$ being fixed, there exists a unique $b=q_b/\phi$ with $q_b$ polynomial in $z$, $\deg_z q_b\leq N+1$ with $b(z)dz=(\frac{-\theta_n(a_i)}{\xi}+o(1))d\xi$ as above and such that the expression $\frac{3}{4}\left(\frac{\partial_z c}{c}\right)^2-\frac{\partial^2_z c}{2c}-\frac{\partial_z b}{2}+\frac{b}{2}\frac{\partial_z c}{c}$ has residues $\nu_i$ at $z=\lambda_i$, $i=1,\ldots,N$. Determination of this $b$ can be explicitly performed by Lagrange interpolation.
 
For these determined $b$ and $c$, there exists a unique $a=q_a/\phi$ with $q_a$ polynomial in $z$, $\deg_z q_a\leq N+1$ such that the corresponding expression $(\ref{p1})$ has principal term $\frac{a_i}{(z-t_i)^2}$ at $z=t_i$. Again notice that, for applications, its expression could be easily exhibited.

With such choices of $a,b,c$, by formula $(\ref{p1})$ we determine an element $p$ that has the same coefficients as $p_0$ in the basis of Lemma $\ref{linalg}$, both must be equal. \emph{The family $(\nabla_t)$ given by these $a,b,c$ will be denoted $\left(\nabla_t(p_0,\theta_n)\right)$}.
~\\

 Moreover, fixing $t$ in the expression of $p$, the pole $\lambda_i$ for  $\tilde{Z}\mapsto \left [\begin{smallmatrix} 0&1\\ p(z,t)&0\end{smallmatrix}\right] dz\cdot \tilde{Z}$
 is apparent if and only if
 \begin{equation}\label{ham}L_i=M_i\left(\sum_{k=1}^N M^{k,i}\nu_k^2- M^{k,i,0}\nu_k -M^{k,i}U_k\right),~~~~i=1,\ldots,N ;
 \end{equation}
 where, for $k,i\in\{1,\ldots,N\}$, using formulae $(\ref{defphipsi})$ for $\psi_{\lambda}$ and $\phi_t$,
 \[\left [ \begin{array}{l}
 M_i=-\frac{\psi_{\lambda}(t_i)}{(\partial_z\phi_t)(t_i)}~;\\
 \\
 M^{k,i}=\frac{\phi_t(\lambda_k)}{(\lambda_k-t_i)(\partial_z\psi_{\lambda})(\lambda_k)}~;\\
 \\
 M^{k,i,0}=M^{k,i} \left ( \sum\limits_{m=1, m\neq i}^{N+2}\frac{1}{\lambda_k-t_m}-\sum\limits_{m=1,m\neq k}^{N}\frac{1}{\lambda_k-\lambda_m}\right );\\
 \\
 U_k=\frac{a_{N+3}}{\lambda_k(\lambda_k-1)}+\sum\limits_{m=1,m\neq i}^{N+2}\frac{a_m}{(\lambda_k-t_m)^2}+\sum\limits_{m=1,m\neq k}^{N}\frac{3}{4(\lambda_k-\lambda_m)^2}~\cdot
 \end{array} \right.\]
 For this equivalence, compare \cite[Prop. $4.3.4$]{MR1118604}.

 \begin{theorem}\label{isomGarnier}
 Fix $N>0$ and let $T$ be a simply connected neighborhood of $t^0$ in $F_{3,N}\Pu$.
 Fix a family of complex numbers $(a_i)_{i=1,\ldots, N+3}$. Fix holomorphic maps $\lambda_i, \nu_i,L_i: T \rightarrow \C$, $i=1,\ldots,N$.
 Define $p(z,t)$ by equation $(\ref{eqp})$.
Consider the family of matrices $P_t=\left [\begin{smallmatrix} 0&1\\ p(z,t)&0\end{smallmatrix}\right]$.
Suppose, for every $t\in T$, the connection $Z\mapsto P_t dz\cdot Z$  has $N+3$ distinct non-apparent poles in $z=t_1,\ldots,t_N,0,1,\infty$ and $N$ distinct apparent poles in $z=\lambda_i(t)$, $i=1,\ldots,N$.

Then $(L_i)$ satisfies  formulae $(\ref{ham})$ and the following are equivalent.
\begin{enumerate}
\item \label{isomfuchs}There exists a matrix  $\Omega=\sum \Omega_i(z,t)dt_i$ of meromorphic $1$-forms on $T\times \Pu$ such that $\omega:=Pdz+\Omega$ satisfies $d\omega=\omega\wedge\omega$. 
\item The functions $\lambda_i(t),\nu_i(t)$ are holomorphic and satisfy the Hamiltonian system
\begin{equation}\label{Garniersystem}
\left \lbrace
 \begin{array}{l}
 \partial_{t_i}\lambda_k= \partial_{\nu_k}L_i,\\
 \partial_{t_i}\nu_k= -\partial_{\lambda_k}L_i, ~~~~~~k,i=1,\ldots,N;
\end{array}
\right .
\end{equation}
\end{enumerate}
 \end{theorem}
 The system $(\ref{Garniersystem})$ is called the Garnier system.\\
{\it On the proof of Theorem $\ref{isomGarnier}$} This theorem was first discovered by Garnier in \cite[Troisi\`eme Partie]{MR1509146}.
 The Hamiltonian formulation and a more detailed proof were given by K. Okamoto in \cite{MR866050} and subsequently in \cite{MR1118604}, to which we refer as the most detailed exposition on the topic.
 
  Yet, in the version of Theorem $\ref{isomGarnier}$ given in \cite{MR866050} and \cite{MR1118604}, our assumption of non apparence for the poles $t_i$ is replaced by the stronger assumption $\theta_i\not \in \Z$, $i=1,\ldots,n$. The only gap in \cite{MR1118604} for a complete proof of Theorem $\ref{isomGarnier}$ is filled by a slight generalization of \cite[Lemma $4.4.2$]{MR1118604} that we will present in the Appendix.
 \qed
 \begin{remark}
 Notice that the first line of $(\ref{Garniersystem})$ allows to derive $\nu_k=(\partial_{t_i}\lambda_k+M_iM^{k,i,0})/(2M_iM^{k,i})$, we can thus eliminate the variable $\nu_k$ from $(\ref{Garniersystem})$ and obtain a system of equations for the sole $N$-tuple $(\lambda_i(t))$.
 This gives the original form of the Garnier system.
 \end{remark}
 \begin{remark}
 Actually, in Theorem $\ref{isomGarnier}$, the coefficients of the entries of $\Omega$ (and obviously of $P(z,t)$) are determined as rational functions of  $(t_i), (\lambda_i(t))$ and their derivatives.  
 However, it is not clear that for every solution $(\lambda_i)$ of the Garnier system, the so determined family $(Z\mapsto dZ-P(z,t)dz\cdot Z)_t$ has a non-apparent pole at every $t_i$. \end{remark}
 
 \begin{definition}
 Let $\left(P(z,t)\right)_t=\left(\left [\begin{smallmatrix} 0&1\\ p(z,t)&0\end{smallmatrix}\right]\right)_t$ be a family which satisfies condition $(\ref{isomfuchs})$ of Theorem $\ref{isomGarnier}$ and denote $(\lambda_i)_i={1,\ldots,N}$ the corresponding solution of the Garnier system.
 We will say that $(\lambda_i)$ governs the isomonodromic deformation $(\nabla_t)$, for $(\nabla_t)=\left(\nabla_t(p,\theta_n)\right)$.
 \end{definition}

\begin{COR:2}
Let $(\lambda_i)$ be a  solution of a Garnier system governing the isomonodromic deformation of a rank $2$ trace free logarithmic connection $\nabla_{t^0}$ on $\Pu$ with no apparent pole.
The following are equivalent.
\begin{enumerate}
\item \label{appli1}The multivalued functions $\lambda_i$  are algebraic functions.
\item \label{appli2}The functions $\lambda_i$ have finitely many branches.
\item \label{appli3} The conjugacy class $[\rho]$ of the monodromy representation $\rho$ of $\nabla_{t^0}$ has finite orbit under $\mathrm{MCG}_n\Pu$.
\end{enumerate}
\end{COR:2}
\begin{proofbis}
The implication $\ref{appli1}.\Rightarrow \ref{appli2}.$ is obvious.
~\\
\noindent $\ref{appli2}.\Rightarrow\ref{appli3}$. The finite branching solution $(\lambda_i)$ is meromorphic on a finite branched cover $\chi:U \rightarrow F_{3,N}\Pu$,
Let $b\in U$ be such that $\chi(b)=t_0$.

Up to changing the base point $t_0$, we may suppose $\chi$ is \'etale  in the neighborhood of $t_0$ and that there exists an (analytic) Zariski neighborhood  $V$ of $b$ in $U$ such that the restriction to $V$ of the function $(\lambda_1,\ldots,\lambda_N,t_1,\ldots,t_N,0,1,\infty)$ is holomorphic and takes values in $F_{3,2N}\Pu$.
By continuation of the local isomonodromic family $\nabla_t(p)$ determined in the neighborhood of $b$ by $(\lambda_i)$, we get a flat connection $\nabla$ on $V\times \Pu$ such that  $\nabla_{\vert \{b\}\times \Pu}$ is isomorphic to $\nabla_{t^0}$ (up to the obvious identification  $\{b\}\times \Pu\simeq \{t^0\}\times \Pu$).
Analogously to Corollary $\ref{finiteindex}$, as $\chi_* \pi_1(U)$ is a finite index subgroup of $\pi_1(F_{3,N}\Pu)$ (\cite[Lemma 4.19]{MR1841091}), the conjugacy class $[\rho]$ has finite orbit under $\pi_1(F_{3,N}\Pu)$. This gives the conclusion, because $\pi_1(F_{3,N}\Pu)$ has index two in $\pi_1(F_{0,n}\Pu)$.
~\\
\noindent $\ref{appli3}.\Rightarrow \ref{appli1}$. It is easily seen, thanks to Theorem $\ref{thmPDL}$ and Lemma $\ref{original}$, that any rank two flat connection with no apparent pole has only mild transversal models.  In particular, we may apply Theorem $\ref{maintheorem}$ to $\nabla_{t^0}$ and infer that the isomonodromic deformation of $\nabla_{t^0}$ is algebraizable. We conclude by Theorem $\ref{AlgGD}$.
\qed \end{proofbis}

\section{Fields of definition}\label{Fields}
The implication $(\ref{unthm})\Rightarrow (\ref{deuxthm})$ of Theorem $\ref{maintheorem}$ can be interpreted as a method to construct some logarithmic connections on ruled varieties. One of the simplest invariants we can retain from such a connection is the set of fields of definition of its projective monodromy group.
Proposition $\ref{proptraces}$ and Corollary $\ref{cortraces}$ below should be helpful to understand these fields in our construction.

\begin{definition}
Let $k$ be a field.
Let $\rho : \Gamma \rightarrow \mathrm{PGL}_m(k)$ be a group morphism and $f$ be a subfield of $k$.
Let $G_f$ be the image of the subgroup $\mathrm{GL}_m(f)<\mathrm{GL}_m(k)$ by $\mathrm{P} : \mathrm{GL}_m(k) \rightarrow \mathrm{PGL}_m(k)$.
 If there exists $C\in \mathrm{PGL}_m(k)$ such that, $C\rho(\Gamma) C^{-1}< G_f$, we say that $f$ is a \emph{field of definition} for $\rho$.
\end{definition}
\begin{proposition}$\label{proptraces}$
Let $\Gamma=N\rtimes Q$ be a semidirect product.
Let $k$ be an algebraically closed field and let $f$ be a subfield of $k$.
Suppose we have a morphism $\rho : \Gamma \rightarrow \mathrm{GL}_m(k)$ such that
 the restriction $\rho_{\vert N}$ to $N$ is irreducible and $f$  is a field of definition for $\mathrm{P}\rho_{\vert N}$.

If there exists a generating set $S\subset N$ for $N$ such that, for every $\alpha\in S$, $\rho(\alpha)$ has nontrivial trace,
   then $f$ is a field of definition for $\mathrm{P}\rho$.
  \begin{proofbis} 
After conjugation, we may suppose   $\mathrm{P}\rho_{\vert N}$ takes values in $G_f$. Fix $\beta \in Q$.

  For $\alpha \in S$, Let $M_{\alpha}$ be an element of $\mathrm{GL}_m(f)$ with $\mathrm{P}M_{\alpha}=\mathrm{P}\rho(\alpha)$. Let $Q_{\alpha}$ be an element of $\mathrm{GL}_m(f)$ with $\mathrm{P}Q_{\alpha}=\mathrm{P}\rho(\beta^{-1} \alpha \beta)$ and $\mathrm{trace}(Q_{\alpha})= \mathrm{trace}(M_{\alpha})$.
By Lemma $\ref{lemsemidirect}$, the element $\mathrm{P}\rho(\beta)$ is characterized by the system of equations $$  \mathrm{P}\rho_{\vert N}(\beta^{-1} \alpha \beta)=\mathrm{P}\rho(\beta)^{-1}\mathrm{P}\rho_{\vert N}(\alpha)\mathrm{P}\rho(\beta),~~\alpha \in S.$$

The matrices $M_{\alpha}$ having nontrivial trace, this implies that the solutions $M\in \mathrm{Mat}_m(k)$ of the system $$M Q_{\alpha}=M_{\alpha}M,~~\alpha \in S$$
are exactly the matrices of the form $M=\lambda \rho(\beta),~\lambda \in k$.

This system can be interpreted as a homogeneous linear system in the entries of $M$, with coefficients in $f$. Therefore, the existence of the nontrivial solution $\rho(\beta)$ ensures the existence of a nontrivial solution with coefficients in $f$, a solution $M\in \mathrm{GL}_m(f)$. This implies $\mathrm{P}\rho(\beta)\in G_f$.
\qed \end{proofbis}
\end{proposition}
\begin{corollary}$\label{cortraces}$
Let $\Gamma=N\rtimes Q$ be a semidirect product, with $N$ finitely generated.
Let $k$ be an algebraically closed field and let $f$ be a subfield of $k$.
Let $\rho : \Gamma \rightarrow \mathrm{GL}_m(k)$ be a morphism such that
 the restriction $\rho_{\vert N}$ to $N$ has Zariski dense image in $\mathrm{GL}_m(k)$ and $f$  is a field of definition for $\mathrm{P}\rho_{\vert N}$.

   The field $f$ is a field of definition for $\mathrm{P}\rho$.
  \begin{proofbis} First, Zariski density of $\rho_{\vert N}$ implies its irreducibility. We want to exhibit a generating set $S$ for $N$ which satisfies the hypothesis of Proposition $\ref{proptraces}$.

  Let $S_0$ be a finite generating set for $N$. 
  The subset $$U:=\{M\in \mathrm{GL}_m(k) \vert \forall \alpha \in S_0,  \mathrm{trace}(M^{-1} \rho(\alpha))\neq 0 \mbox{ and }  \mathrm{trace}(M)\neq 0 \}$$
  is a Zariski open subset of  $\mathrm{GL}_m(k)$. By density, it must contain an element $\rho(\alpha_0)$, with $\alpha_0\in N$.
  The generating set $S:=\{\alpha_0\} \cup \{\alpha_0^{-1}\alpha, \alpha \in S_0\}$ satisfies the hypothesis of Proposition $\ref{proptraces}$.
  \qed \end{proofbis}
\end{corollary}

\appendix
\section*{Appendix}\label{apGarnier}

We want to explain here the generalization of Lemma $4.4.2$ of \cite{MR1118604} we need to obtain Theorem $\ref{isomGarnier}$ with our slightly weakened hypotheses.
The generalized lemma is as follows.
In the sequel we use the terminology of \cite{MR1118604}, beware that the wording \emph{``logarithmic pole"} refers to the form of the local solution of the considered Fuchsian
scalar equation.

\begin{lem:A1}
Take an integrable system $$dY=\left(\left [\begin{smallmatrix} 0&1\\ p(z,t)&0\end{smallmatrix} \right ]dz+\Omega\right)\cdot Y$$
where $\Omega=\sum\limits_{i=1}^N\Omega_i(z,t) dt_i$ and $p$ is given by equation $(\ref{eqp})$.
Suppose $z\mapsto p(z,t^0)$ has exactly $2N+2$ poles in $\C$ (without multiplicity) and the monodromy of the system around any of the $N+3$ poles $z=0,1,\infty$ and $z=t_i,i=1,\ldots ,N$ is not scalar (non-apparent poles).
Then the $(1,2)$ entry $A_i$ of $\Omega_i$ satisfies the following.

The rational function $z\mapsto A_i(z,t^0)$
\begin{enumerate}
\item  is holomorphic outside $\lambda_1,\ldots,\lambda_N,\infty$,
\item has only simple poles in $\C \cup \{\infty\}$ and
\item has zeroes in the points $0,1$ and $t_j^0$, for $1\leq j \leq N$, $j\neq i$.
\end{enumerate}

\begin{proofbis}
The proof is given in \cite{MR1118604} with the hypothesis that none of the poles $z=0,1,\infty,t_1,\ldots,t_N$ has integer local exponent $\theta$.
It is obtained through local studies in the points $z=0,1,\infty,t_i,\lambda_i(t)$ involving local solutions of $\frac{d^2y}{dz^2}=p(z,t)y$ given by Frobenius's method.

To complete the proof, we only need to do the analogous study at those elements of $\{0,1,\infty,t_1,\ldots,t_N\}$ that are logarithmic poles. 

Suppose $z=t_j$ is a logarithmic pole of $\frac{d^2y}{dz^2}=p(z,t)y$, with $t_j\neq t_i$, $t_j\neq \infty$.

By the arguments of \cite{MR1118604}, equation $4.4.5$ p. $185$, we know 
\begin{eqnarray}\label{eqAi} A_i(z,t)=u(t)g(z,t)+a(t)v_1^2+b(t)v_1v_2+c(t)v_2^2
\end{eqnarray}
for
\begin{itemize}
 \item $v_1(z,t),v_2(z,t)$ any fundamental system of solutions of $\frac{d^2y}{dz^2}=p(z,t)y$ near $z=t_j$,
\item $g(z,t):=v_1\frac{\partial}{\partial t_i}v_2-v_2\frac{\partial}{\partial t_i}v_1$
\item $u(t), a(t),b(t),c(t)$ holomorphic functions of $t$, with $u$ nowhere vanishing.
\end{itemize}
By Frobenius's method (with parameter) we may take 
\begin{itemize}
\item $v_1=(z-t_j)^{\frac{1}{2}}h_1(z,t),$\\ $v_2=(z-t_j)^{\frac{1}{2}}(h_1(z,t)\log(z-t_j)+h_2(z,t))$, ~~ if $\theta_j=0$;
\item $v_1=(z-t_j)^{\frac{1+m}{2}}h_1(z,t)$,\\ $v_2=(z-t_j)^{\frac{1-m}{2}}h_1(z,t)+k(t)(v_1 \log(z-t_j)+(z-t_j)^{\frac{1-m}{2}}h_3(z,t))$,
~~ if $\vert \theta_j \vert=m>0$.
\end{itemize}
Where the functions $h_1,h_2,h_3$ are holomorphic in the neighborhood of $z=t_j$, $h_1$ does not vanish and $k(t)$ is a nowhere vanishing holomorphic function of $t$.

In any case, the right hand side of $(\ref{eqAi})$
takes the form  $$\sum_{\ell=0}^2 k_{\ell}(z,t)\left(\log (z-t_j)\right)^{\ell},$$ with $k_{\ell}(z,t)$ holomorphic in the neighborhood of $z=t_j$. The lemma below and uniformity of $A_i$ allow to conclude $k_2=k_1=0$. In any case, the specific form of $k_0$ and the conditions given by $k_2=k_1=0$ allow to see $k_0(z,t)$ vanishes at $z=t_j$.

The same arguments allow to prove holomorphicity of $A_i$ at $t_i$, in case it is a logarithmic pole.
If $\infty$ is a logarithmic pole  of $\frac{d^2y}{dz^2}=p(z,t)y$, it remains to see that $A_i$ has at most a pole of order one at this point.
The  proof is almost as above, the only difference is a slight change in the form of the local fundamental system of solutions.
\qed \end{proofbis}
\end{lem:A1}
\begin{lem:A2}
Let $(k_{\ell}(z))_{0\leq\ell\leq r}$ be holomorphic functions in the punctured disc $\D^*:=\{z\in \C^*,\vert z \vert<1\}$.
Consider the function \[f(z)=\sum_{\ell=0}^r k_{\ell}(z)\left(\log z\right)^{\ell}\] defined on the universal cover of $\D^*$.
If $f$ is uniform then $k_{\ell} \equiv 0$, $\ell>0$.
\begin{proofbis}
We proceed by induction on $r$. The result is trivial for $r=0$.
Take $r>0$ and suppose the lemma is true for $r-1$.  The monodromy group of $\log z$ is generated by $\log z \mapsto \log z+2i\pi$. Hence, $f$ is uniform if and only if
the function $g(z)=f(z)-\sum_{\ell=0}^r k_{\ell}(z)\left(\log z +2i\pi \right)^{\ell}$ is zero.
This function $g(z)$ is obviously uniform and can be written in the form $g(z)=\sum_{\ell=0}^{r-1} q_{\ell}(z)\left(\log z \right)^{\ell}$ with $q_{\ell}$ holomorphic in $\D^*$.
By the induction hypothesis $2i\pi rk_r=q_{r-1}\equiv 0$. Again by the induction hypothesis, $k_{\ell}\equiv 0$ for $\ell=1,\ldots, r-1$.

\qed \end{proofbis}
\end{lem:A2}


\bibliographystyle{amsalpha}
\bibliography{biblio}
\end{document}